\documentclass[letterpaper, 12pt]{amsart}
\pdfoutput=1 
\usepackage{MA}

\title
[Boundary symmetries of topological orders]
{Boundary symmetries of (2+1)D topological orders}

\author{Kylan Schatz}

\begin{document}

    \begin{abstract}
        We elaborate an algebraic framework for describing internal topological symmetries of gapped boundaries of (2+1)D topological orders. We present a categorical obstruction to the coherence of bulk group symmetry and boundary symmetries in terms of liftings of categorical actions on the bulk theory to a certain 2-group of boundary symmetries.
    \end{abstract}

    \maketitle

    \section{Introduction}

    \renewcommand{\thetheorem}{\Alph{theorem}}

    The study of symmetry enriched topological order (SETO) is an important component of the field of topological phases of matter, with applications in both condensed matter physics and quantum information \cite{BBCW19, Del19}. Corresponding to a gapped (2+1)D topological phase enriched with a global on-site $G$-symmetry is a canonical SETO described by a $G$-crossed braided extension of a modular tensor category (MTC).
    
    An interesting class of topological phase transitions between gapped (2+1)D phases are those which arise from \emph{anyon condensation} \cite{Kon14, Bur18}. Mathematically, condensable anyons are characterized by commutative algebra objects in the bulk MTC. Understanding the phase transitions which preserve symmetry is a topic of additional import \cite{BJLP19}. Obstructions to the preservation of symmetry are equivalent to lifts of the underlying group action to the 2-group of boundary symmetries. Further, equivariant structure on the corresponding algebra object is known to be sufficient for preservation of group symmetry. 

    The boundary theory has its own internal topological symmetries - namely `string operators' - which are characterized by matrix product operators (MPOs) \cite{BMW+17, WBV17}. These symmetries holographically represent the `dual' string net model, and their algebraic structure describes a \emph{hypergroup} \cite{SW03}. It is natural to ask when extensions of bulk symmetry to the boundary are coherent with internal symmetries of the boundary theory. 
   
    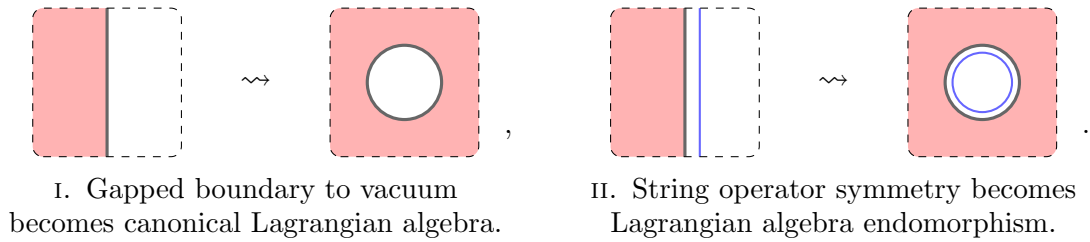
\begin{figure}[ht]
        \centering
        \begin{subfigure}[h]{0.4\textwidth}
            \centering
            \adjustbox{max width=0.9\textwidth}{    \begin{tikzpicture}
        \draw[rounded corners, draw=none, fill=red!30] (1, 0) -- (0, 0) -- (0, 2) -- (1, 2);
        \draw[very thick, draw=black!60] (1, 0) -- (1, 2);
        \draw[rounded corners, dashed] (0, 0) -- (2, 0) -- (2, 2) -- (0, 2) -- cycle;

        \draw[rounded corners, dashed, fill=red!30] (4, 0) -- (6, 0) -- (6, 2) -- (4, 2) -- cycle;
        \draw[anchor=center, very thick, draw=black!60, fill=white] (5, 1) circle (0.5);

        \node at (3, 1) {$\rightsquigarrow$};
    \end{tikzpicture}}
            \caption{Gapped boundary to vacuum becomes canonical Lagrangian algebra.}
        \end{subfigure}, \qquad
        \begin{subfigure}[h]{0.4\textwidth}
            \centering
            \adjustbox{max width=0.9\textwidth}{    \begin{tikzpicture}
        \draw[rounded corners, draw=none, fill=red!30] (1, 0) -- (0, 0) -- (0, 2) -- (1, 2);
        \draw[very thick, draw=black!60] (1, 0) -- (1, 2);
        \draw[thick, draw=blue!60] (1.2, 0) -- (1.2, 2);
        \draw[rounded corners, dashed] (0, 0) -- (2, 0) -- (2, 2) -- (0, 2) -- cycle;

        \draw[rounded corners, dashed, fill=red!30] (4, 0) -- (6, 0) -- (6, 2) -- (4, 2) -- cycle;
        \draw[anchor=center, very thick, draw=black!60, fill=white] (5, 1) circle (0.5);
        \draw[anchor=center, thick, draw=blue!60, fill=white] (5, 1) circle (0.4);

        \node at (3, 1) {$\rightsquigarrow$};
    \end{tikzpicture}}
            \caption{String operator symmetry becomes Lagrangian algebra endomorphism.}
        \end{subfigure}.
        \caption{`Rolling up' the boundary.}
        \label{heuristic}
    \end{figure}

    For a separable algebra $A$ object in a unitary braided fusion category one may endow its endomorphism ring $\End(A)$ with a \emph{convolution product} -- a generalization of the convolution of finite group characters -- to obtain the \emph{convolution algebra}. If $A$ is commutative, its convolution algebra is commutative semisimple with respect to this convolution product \cite{BD18}. From the data of the convolution algebra of a commutative separable algebra $A$ we may thus define a hypergroup called $\hAut(A)$ with extreme points given by orthogonal convolution idempotents of $\End(A)$. There is a natural identification of the hypergroup of string operators on a boundary spatially representing condensation by $A$ with $\hAut(A)$. When $A \in \mc C$ is \emph{Lagrangian}, $\hAut(A)$ can be identified with the fusion ring of $\mc C_A$ in a canonical way \cite{BJ21}.

    Let $\mc F$ a $G$-crossed braided extension of unitary modular tensor category $\mc E$ with braided action $G\curvearrowright \mc E$. Equivariant structure on the algebra object $A$ provides for each $g$-graded defect $X \in \mc F_g$ an isomorphism which is coherent and compatible with multiplication \cite{BJLP19}:
    \begin{align*}
        \adjustbox{scale=0.7}{    \begin{tikzpicture}[baseline=(current bounding box.center)]
        \draw (1, 0) node[label=below:$A$] {}  -- (1, 1) .. controls (1, 1.5) and (0, 1.5) .. (0, 2)  -- (0, 2.5) node[draw=black, fill=white] {$\lambda^{-1}$} -- (0, 3) node[label=above:$A$] {};
        \path[fill=white] (0.5, 1.5) circle (4pt);
        \draw (0, 0) node[label=below:$X$] {} -- (0, 1) .. controls (0, 1.5) and (1, 1.5) .. (1, 2) -- (1, 3) node[label=above:$X$] {};

        \draw[dashed] (-0.5, 1) -- (1.5, 1) -- (1.5, 1.5) node[label=right:$\beta_{X, A}$] {} -- (1.5, 2) -- (-0.5, 2) -- cycle;
    \end{tikzpicture}} = (\lambda^{-1} \otimes \id) \circ \beta_{X, A}.
    \end{align*}

    String operator symmetries should be compatible with these isomorphisms:
    \begin{align*}
        \adjustbox{scale=0.7}{    \begin{tikzpicture}[baseline=(current bounding box.center)]
        \draw (0, 0) node[label=below:$X$] {} -- (0, 1) .. controls (0, 1.5) and (1, 1.5) .. (1, 2) -- (1, 3) node[label=above:$X$] {};
        \path[fill=white] (0.5, 1.5) circle (4pt);
        \draw (1, 0) node[label=below:$A$] {}  -- (1, 0.5) node[draw=black, fill=white] {$\omega$} -- (1, 1) .. controls (1, 1.5) and (0, 1.5) .. (0, 2) -- (0, 2.5) node[draw=black, fill=white] {$\lambda^{-1}$}  -- (0, 3) node[label=above:$A$] {};
    \end{tikzpicture}} = \,\,
        \adjustbox{scale=0.7}{    \begin{tikzpicture}[baseline=(current bounding box.center)]
        \draw (0, 0) node[label=below:$X$] {} -- (0, 1) .. controls (0, 1.5) and (1, 1.5) .. (1, 2) -- (1, 3) node[label=above:$X$] {};
        \path[fill=white] (0.5, 1.5) circle (4pt);
        \draw (1, 0) node[label=below:$A$] {} -- (1, 1) .. controls (1, 1.5) and (0, 1.5) .. (0, 2) -- (0, 2.1) node[draw=black, fill=white] {$\lambda^{-1}$} -- (0, 2.7) node[draw=black, fill=white] {$\omega$}  -- (0, 3) node[label=above:$A$] {};
    \end{tikzpicture}},
    \end{align*}

    which is equivalent to the following diagram commuting:
    \begin{center}
        \begin{tikzcd}
            A 
                \arrow[r, "\omega"]
                \arrow[d, "\lambda"]
            & A
                \arrow[d, "\lambda"] \\
            \alpha(A)
                \arrow[r, "\alpha(\omega)"]
            & \alpha(A)
        \end{tikzcd}.
    \end{center}

    Further, if a braided action does not admit equivariant structure on the algebra $A$ which is compatible with string operator symmetries, then internal topological symmetries are incompatible with the phase transition represented by the algebra $A$ in any microscopic realization of the associated SETO. Given that the string operator symmetries are represented by hypergroup $\sf H$, we may construct a categorical group representing compatible autoequivalences:
    \begin{align*}
        \underline{\Aut_\otimes^{br}}(\mc E \vert A, \sf{H}) 
            & := \s{(\alpha, \eta^\alpha, \lambda^\alpha) \in \underline{\Aut_\otimes^{br}}(\mc E \vert A) : \lambda^\alpha \circ e_i = \alpha(e_i) \circ \lambda^\alpha, \quad \forall e_i \in \sf{H}}.
    \end{align*}

    Additionally, given a unitary fusion category $\mc C$ and a full unitary fusion subcategory $\mc D \subseteq \mc C$, we may construct a categorical group of `trivializable' autoequivalences:
    \begin{align*}
        \underline{\Aut_\otimes}(\mc C \vert_{\mc D}) 
            & := \s{(\alpha, \eta^\alpha) \in \underline{\Aut_\otimes}(\mc D) : \alpha\vert_{\mc D} \cong \Id_{\mc D} \text{ as \emph{linear} functors}}.
    \end{align*}
    
    \begin{theorem}
        Let $\mc C$ unitary fusion, $L = I(\mathbb 1)$ its canonical Lagrangian algebra, and $\mc D \subseteq \mc C$ full fusion. Then, there is an equivalence of 2-groups:
        \begin{align*}
            \underline{I^L_{\mc D}}:\underline{\Aut_\otimes}(\mc C \vert_{\mc D}) \xrightarrow{\sim} \underline{\Aut_\otimes^{br}}(\mc Z(\mc C) \vert L, \sf{K_0}(\mc D)).
        \end{align*}
    \end{theorem}



    \section{Preliminaries}

    \renewcommand{\thetheorem}{\arabic{section}.\arabic{theorem}}

    \subsection{Monoidal Categories}

    For a category $\mc C$, we will denote the class of objects of $\mc C$ by $\Ob (\mc C)$ and the morphisms from $X$ to $Y$ by $\Hom_{\mc C}(X, Y)$. We denote identity morphisms $\id_X \in \Hom_{\mc C}(X, X)$ and the identity functor $\Id_{\mc C}:\mc C \xrightarrow{} \mc C$, sometimes omitting subscripts. For convenience, we assume that all monoidal categories are strict.
    
    Recall that a \emph{monoidal functor} \cite[Section~2.4]{EGNO} is a tuple $(\alpha, \eta^\alpha)$ of a functor $\alpha$ between monoidal categories and a natural isomorphism:
     \begin{align} \label{mdist}
         \eta^\alpha:\alpha(\blank) \otimes \alpha(\blank) \xrightarrow{} \alpha(\blank \otimes \blank),
     \end{align}
    such that the following diagram commutes:
    \begin{equation} \label{massc}
        
    \begin{tikzcd}
        & \alpha(X) \otimes \alpha(Y) \otimes \alpha(Z) 
            \arrow[dl, "\eta^\alpha_{X, Y} \otimes \id_Z"]
            \arrow[dr, swap, "\id_X \otimes \eta^\alpha_{Y, Z}"]
            \\
        \alpha(X \otimes Y) \otimes \alpha (Z)
            \arrow[dr, "\eta^\alpha_{X \otimes Y, Z}"]
        && \alpha(X) \otimes \alpha(Y \otimes Z)
            \arrow[dl, swap, "\eta^\alpha_{X, Y \otimes Z}"]\\
        & \alpha (X \otimes Y \otimes Z)
    \end{tikzcd}
.
    \end{equation}

    A \emph{monoidal natural transformation} $\pi:(\alpha, \eta^\alpha) \xrightarrow{} (\beta, \eta^\beta)$ is a natural transformation $\pi:\alpha \xrightarrow{} \beta$ such that the following diagram commutes:
    \begin{equation*} 
    \label{mnat}
        
    \begin{tikzcd}
        \alpha(X) \otimes \alpha(Y)
            \arrow[r, "\eta^\alpha_{X, Y}"]
            \arrow[d, "\pi_X \otimes \pi_Y"]
        & \alpha(X \otimes Y)
            \arrow[d, "\pi_{X \otimes Y}"]\\
        \beta(X) \otimes \beta(Y) 
            \arrow[r, "\eta^\beta_{X, Y}"]
        & \beta(X \otimes Y)
    \end{tikzcd}
.
    \end{equation*}

    We will refer to naturality in the sense of equation \eqref{mdist} as \emph{distributive} and naturality in the sense of the diagram \eqref{massc} as \emph{associative}.

    A \emph{unitary fusion category} \cite{Pen18} is a finitely semisimple $C^*$-tensor category. We will denote the object dual to $X$ by $\bar X$ and the dagger structure by $\dagger:\Hom(\blank_{(1)}, \blank_{(2)}) \xrightarrow{\sim} \Hom(\blank_{(2)}, \blank_{(1)})$. Unless otherwise stated, we assume that all isomorphisms in a unitary fusion category are unitary. A \emph{unitary monoidal functor} is a monoidal functor $(\alpha, \eta^\alpha)$ between unitary fusion categories for which $\alpha$ is a $\dagger$-functor and $\eta^\alpha$ is a unitary natural isomorphism.

    For $\mc C$ unitary fusion, the \emph{Drinfeld center of $\mc C$} \cite[Definition~3.4]{Mug03} denoted $\mc Z(\mc C)$ has objects $(X, \psi)$ where $X \in \Ob \mc C$ and $\psi:X\otimes\blank \xrightarrow{} \blank\otimes X$ is an isomorphism such that for any $Y, Z \in \Ob \mc C$ and $f \in \Hom_{\mc C}(Y, Z)$:
    \begin{align}
            & \label{br-nat} \psi^\alpha_Z \circ (\id \otimes f) = (f \otimes \id) \circ \psi^\alpha_Y, \quad \text{and} \\
            & \label{br-hex} \psi^\alpha_{Y \otimes Z} = (\id \otimes \psi^\alpha_Z) \circ (\psi^\alpha_Y \otimes \id).
        \end{align}
    In other words, $\psi$ is a natural isomorphism satisfying the braid relation \eqref{br-hex}. Morphisms in the center are compatible with half-braidings, that is:
    \begin{align}
        \label{zc-homs}
        \Hom_{\mc Z(\mc C)}((X, \psi), (Y, \varphi) = \s{f \in \Hom_{\mc C}(X, Y) : (\id \otimes f) \circ \psi = \varphi \circ (f \otimes \id)}.
    \end{align}
    
    The Drinfeld center is monoidal with product:
    \begin{align*}
        (X, \psi) \otimes (Y, \varphi) = (X \otimes Y, (\psi \otimes \id) \circ (\id \otimes \varphi)).
    \end{align*}
    
    The Drinfeld center of a unitary fusion category is unitary modular \cite{Mug03, BV13}. There is a natural monoidal forgetful functor $\Forg:\mc Z(\mc C) \xrightarrow{} \mc C$, denote its (left) adjoint $I:\mc C \xrightarrow{} \mc Z(\mc C)$ the \emph{induction functor}. We will use the following concrete model \cite{BK10, BJ21}, described in the (bottom-up/optimistic) graphical calculus \cite{Coe10}:
    \begin{gather*}
        I(X) = \s[p]{\bigoplus_{U \in \Irr \mc C} U \otimes X \otimes \bar U, \, \psi_{I(X)}}, \qquad I(f) = \bigoplus_{U \in \Irr \mc C} \id_U \otimes f \otimes \id_{\bar U}, \\
        \psi_{I(X), W} = \bigoplus_{U, V \in \Irr \mc C} \sum_i \sqrt{\dim U} \sqrt{\dim V} \; 
            \adjustbox{scale=0.7}{    \begin{tikzpicture}[baseline=(current bounding box.center)]
        \draw (0,0) node[label={below:{\strut$U$}}] {} -- (0, 2) node[label={above:{$V$}}] {};
        \draw (0.5,0) node[label={below:{\strut$X$}}] {} -- (0.5, 2) node[label={above:{$X$}}] {};
        \draw (1,0) node[label={below:{\strut$\bar U$}}] {} -- (1, 2) node[label={above:{$\bar V$}}] {};
        \draw (0, 1) .. controls (-0.5, 1.7) .. (-0.5, 2) node[label={above:{$W$}}] {};
        \draw (1, 1) .. controls (1.5, 0.3) .. (1.5, 0) node[label={below:{\strut$W$}}] {};
        \draw[black, fill=white] (0, 1) circle (2pt) node[label={left:{$i^*$}}] {};
        \draw[black, fill=white] (1, 1) circle (2pt) node[label={right:{$i$}}] {}; 
    \end{tikzpicture}},
    \end{gather*}
    where direct sums are over some choice of representatives of isomorphism classes of simple objects and $\s{i^*} \subseteq \Hom(U, W \otimes V)$, $\s{i} \subseteq \Hom(\bar U \otimes W, \bar V)$ are dual bases with respect to the pairing:
    \begin{align} \label{pepsi}
        \delta_{i, j} = 
        \adjustbox{scale=0.7}{    \begin{tikzpicture}[baseline=(current bounding box.center)]
        \draw (0, 1.5) arc (180:0:0.5);
        \draw (0, 0.5) arc (180:360:0.5);
        \draw (0,0.5) -- (0, 1.5) .. controls (0, 1) and (1, 1) .. (1, 0.5) -- (1, 1.5);
        \draw (0, 1) node[label={west:\strut$\bar U$}] {};
        \draw (1, 1) node[label={east:\strut$V$}] {};
        \draw (0.5, 1) node[label={south:$W$}] {};
        \draw[black, fill=white] (0, 1.5) circle (2pt) node[label={left:{$i$}}] {};
        \draw[black, fill=white] (1, 0.5) circle (2pt) node[label={right:{$j^*$}}] {}; 
    \end{tikzpicture}},
    \end{align}
    using the canonical unitary spherical structure.

    \subsection{Algebra Objects}

    A \emph{Q-system object} \cite[Definition~3.1]{CHJP22} in a unitary fusion category $\mc C$ is a separable $C^*$-Frobenius algebra $(A, \Delta, \iota, \nabla, \epsilon)$ whose comultiplication and counit are the dagger of the multiplication and unit respectively. The \emph{canonical Lagrangian algebra} ${L = I(\mathbb 1)}$ has the structure of a commutative Q-system object in the Drinfeld center $\mc Z(\mc C)$ with structure maps:
    \begin{align*}
        & \Delta_L = \bigoplus_{X} \frac{1}{\sqrt{\dim X}} \,
            \adjustbox{scale=0.7}{    \begin{tikzpicture}[baseline=(current bounding box.center)]
        \draw (0, 0) node[label={below:{\strut$X$}}] {} .. controls (-0.1, 1.75) and (0.75, 1) .. (0.75, 2) node[label={above:{\strut$X$}}] {};
        \draw (2, 0) node[label={below:{\strut$\bar X$}}] {} .. controls (2.1, 1.75) and (1.25, 1) .. (1.25, 2) node[label={above:{\strut$\bar X$}}] {};
        \draw (0.5, 0) node[label={below:{\strut$\bar X$}}] {} -- (0.5, 0.5);
        \draw (1.5, 0) node[label={below:{\strut$X$}}] {} -- (1.5, 0.5);
        \draw (0.5, 0.5) arc (180:0:0.5);
    \end{tikzpicture}},
        & \nabla_L = \Delta_L^\dagger = \bigoplus_{X} \frac{1}{\sqrt{\dim X}} \,
            \adjustbox{scale=0.7}{    \begin{tikzpicture}[baseline=(current bounding box.center)]
        \draw (0, 0) node[label={above:{\strut$X$}}] {} .. controls (-0.1, -1.75) and (0.75, -1) .. (0.75, -2) node[label={below:{\strut$X$}}] {};
        \draw (2, 0) node[label={above:{\strut$\bar X$}}] {} .. controls (2.1, -1.75) and (1.25, -1) .. (1.25, -2) node[label={below:{\strut$\bar X$}}] {};
        \draw (0.5, 0) node[label={above:{\strut$\bar X$}}] {} -- (0.5, -0.5);
        \draw (1.5, 0) node[label={above:{\strut$X$}}] {} -- (1.5, -0.5);
        \draw (0.5, -0.5) arc (180:360:0.5);
    \end{tikzpicture}}, \\
        & \iota_L = \bigoplus_{X} \sqrt{\dim X} \,
            \adjustbox{scale=0.7}{
    \begin{tikzpicture}[baseline=(current bounding box.center)]
        \draw (1, 0) -- (1, 0.5) node[label={above:{\strut$\bar X$}}] {};
        \draw (0, 0) -- (0, 0.5) node[label={above:{\strut$X$}}] {};
        \draw (0, 0) arc (180:360:0.5);
    \end{tikzpicture}
}, 
        & \epsilon_L = \iota_L^\dagger = \bigoplus_{X} \sqrt{\dim X} \, 
            \adjustbox{scale=0.7}{
    \begin{tikzpicture}[baseline=(current bounding box.center)]
        \draw (1, 0) -- (1, -0.5) node[label={below:{\strut$\bar X$}}] {};
        \draw (0, 0) -- (0, -0.5) node[label={below:{\strut$X$}}] {};
        \draw (0, 0) arc (180:0:0.5);
    \end{tikzpicture}
}.
    \end{align*}
    
    Let $A \in \Ob \mc C$ a Q-system object. The \emph{category of right $A$-modules over $\mc C$} $\mc C_A$ {\cite[Section~7.8]{EGNO}} has objects $(X, \mu)$ where $X \in \Ob \mc C$ and $\mu:X \otimes A \xrightarrow{} X$ is a morphism such that the following diagram commutes:
    \begin{equation} \label{mod-action}
        
    \begin{tikzcd}
        & X \otimes A \otimes A
            \arrow[dl, "\id \otimes \Delta"]
            \arrow[dr, swap, "\mu \otimes \id"]\\
        X \otimes A 
            \arrow[dr, "\mu"]
        && X \otimes A
            \arrow[dl, swap, "\mu"]\\
        & X
    \end{tikzcd}
.
    \end{equation}
    
    Morphisms in $\mc C_A$ are compatible with the module action, that is:
    \begin{align} \label{mod-mor}
        \Hom_{\mc C_A}((X, \mu), (Y, \nu)) = \s{f \in \Hom_{\mc C}(X, Y) : f \circ \mu = \nu \circ (f \otimes \id)}.
    \end{align}
    
    When $A$ is commutative and central, $\mc C_A$ has the structure of a unitary fusion category \cite[Proposition~5.1]{BN10}. When $L = I(\mathbb 1)$ the canonical Lagrangian, $I_L: \mc C \cong \mc Z(\mc C)_L$ where $I_L(X) = (I(X), \eta^I_{X, \mathbb 1})$ is an equivalence of unitary fusion categories \cite[Corollary~5.6]{BN10}.

    \subsection{Categorical Groups}

    A \emph{categorical group (2-group)} \cite[Definition~2.11.4]{EGNO} is a 2-category with one object whose morphisms are invertible up to 2-morphism and whose 2-morphisms are invertible. We denote categorical groups and functors of categorical groups with \underline{underline}. A \emph{2-subgroup} is the one object together with morphisms given by a full monoidal subcategory of the original one morphism space.

    Let $\mc C$ a unitary fusion category. $\underline{\Aut_{\otimes}}(\mc C)$ is the categorical group with one object whose morphisms are unitary monoidal autoequivalences of $\mc C$ and whose 2-morphisms are unitary monoidal natural isomorphisms. 

    \begin{definition}[{\cite[Definition~4.1]{BJLP19}}]
        Let $\mc C$ unitary modular category and $A \in \mc C$ a commutative Q-system object. Define $\underline{\Aut_\otimes^{br}}(\mc C \vert A)$ to be the categorical group whose morphisms are tuples $(\alpha, \eta^\alpha, \lambda^\alpha)$, where $(\alpha, \eta^\alpha)$ is a unitary braided monoidal autoequivalence of $\mc C$ and $ \lambda^\alpha : A \xrightarrow{}  \alpha(A)$ an algebra isomorphism; and 2-morphisms are:
        \begin{align*}
            \Hom((\alpha, \eta^\alpha, \lambda^\alpha), (\beta, \eta^\beta, \lambda^\beta)) = \s{\pi \in \Hom((\alpha, \eta^\alpha), (\beta, \eta^\beta)) : \pi_A \circ \lambda^\alpha = \lambda^\beta}.
        \end{align*}
    \end{definition}

    We denote by $\Aut_\otimes(\mc C)$ and $\Aut_\otimes^{br}(\mc C \vert A)$ the group truncations of $\underline{\Aut_\otimes}(\mc C)$ and $\underline{\Aut_\otimes^{br}}(\mc C \vert A)$.

    \begin{remark}
        Let $(F, J):\mc C \xrightarrow{} \mc D$ a unitary monoidal equivalence of unitary fusion categories. Then, $(F, J)$ extends to a canonical unitary monoidal equivalence of 2-groups:
        \begin{align*}
            & \underline{F}:\underline{\Aut_\otimes}(\mc C) \xrightarrow{} \underline{\Aut_\otimes}(\mc D), 
            \qquad \underline{F}(\alpha, \eta^\alpha) := (\tilde\alpha, \tilde\eta^\alpha), \\
            \text{where } \quad & \tilde\alpha \circ F(X) := F \circ \alpha(X),
            \qquad \tilde\alpha\circ F(f) := F \circ \alpha(f), \\
            \text{and } \quad & \tilde\eta^\alpha_{FX, FY} := \tilde\alpha(J_{X, Y})^{-1} \circ F(\eta^\alpha_{X, Y}) \circ J_{\alpha X, \alpha Y},
        \end{align*}
        which acts on 2-morphisms $\pi \in \Hom((\alpha, \eta^\alpha), (\beta, \eta^\beta))$ by:
        \begin{align*}
            \underline{F}(\pi)_{FX} = F(\pi_X) .
        \end{align*}
    \end{remark}

    \section{The 2-Group of Lagrangian-Algebra Equivariant Autoequivalences}

    Well known to experts is the correspondence between modular tensor categories containing a Lagrangian algebra object and fusion categories. There is an extension of this correspondence to the level of 2-groupoids \cite{DMNO13}: 1-morphisms on pairs of MTCs and Lagrangian are algebras carry an algebra isomorphism fixing the Lagrangian algebra which is compatible with 2-morphisms. 

    Suppose $\mc C$ is a unitary fusion category of rank $n$ and $(\alpha, \eta^\alpha)$ is a monoidal autoequivalence of $\mc C$. Then, $\alpha$ need descend to a set permutation of simple objects $\sigma \in S([\Irr \mc C]) \cong S(n)$. As $\mc C$ is an Abelian category, there is a natural isomorphism:
    \begin{align}
        \label{irrperm}
        & S:F_{\id} \xrightarrow{} F_{\sigma}, \qquad \text{where} \qquad F_{\id}, F_{\sigma}:\mc C^{\times n} \xrightarrow{} \mc C^{\oplus n} \qquad \text{by} \\ 
        \notag
        & F_{\id}(X_1, X_2, \dots, X_n) = X_1 \oplus X_2 \oplus \dots \oplus X_n, \\ 
        \notag
        & F_{\sigma}(X_1, X_2, \dots, X_n) = X_{\sigma(1)} \oplus X_{\sigma(2)} \oplus \dots \oplus X_{\sigma(n)}.
    \end{align}

    Choose component isomorphisms $\lambda^\alpha_X:Y \xrightarrow{\sim} \alpha(X)$ such that $(\lambda^\alpha_{\bar X})^{-1} = (\lambda^\alpha_{X})^{*}$. Define map:
    \begin{align}
        \label{algisom}
        \tilde\lambda^\alpha = \s[p]{\bigoplus_{U \in \Irr \mc C} \eta^{\alpha}_{U, }} \circ \s[p]{\bigoplus_{U \in \Irr \mc C} \lambda^\alpha_U \otimes \lambda^\alpha_{\bar U}} \circ S_L.
    \end{align}

    So, there is a canonical isomorphism $L \xrightarrow{} \alpha(L)$ making $\alpha(L)$ into an algebra. Choose unitary monoidal natural isomorphism $\xi^\alpha:\Id_{\mc C} \xrightarrow{\sim} \alpha\alpha^{-1}$. Given object $(X, \psi) \in \mc Z(\mc C)$, define a collection of maps $\psi^\alpha = \s{\psi^\alpha_Y:\alpha (X) \otimes Y \xrightarrow{\sim} Y \otimes \alpha (X)}$ by:
    \begin{align}
        \label{newbraid}
        \psi^\alpha_Y = ((\xi^\alpha_Y)^{-1}\otimes \id_{\alpha X}) \circ (\eta^\alpha_{\alpha^{-1}Y, X})^{-1} \circ \alpha(\psi_{\alpha^{-1}Y}) \circ \eta^\alpha_{X, \alpha^{-1}Y} \circ (\id_{\alpha X} \otimes \xi^\alpha_Y).
    \end{align}

    \begin{lemma}
        Define a map $\underline{I^L}:\underline{\Aut_\otimes}(\mc C) \xrightarrow{} \underline{\Aut_\otimes^{br}}(\mc Z(\mc C) \vert L)$ such that $\underline{I^L}(\alpha, \eta^\alpha) := (\tilde\alpha, \tilde\eta^\alpha, \tilde\lambda^\alpha)$ and $\underline{I^L}(\pi)_X := \pi_{\Forg X}$ where $\tilde\alpha(X, \psi) = (\alpha X, \psi^\alpha)$ as in \eqref{newbraid}, $\tilde\alpha(f) = \alpha \circ \Forg(f)$, $\tilde\eta^\alpha = \eta^\alpha$, and $\tilde\lambda^\alpha$ as in \eqref{algisom}. The map $\underline{I^L}$ is an equivalence of 2-groups.
    \end{lemma}

    \begin{proof}
        Let $(\alpha, \eta^\alpha)$ a monoidal autoequivalence of $\mc C$. We first show that $\psi^\alpha$ is a half braiding. Consider the diagram in figure \ref{braid-nat}.
        \begin{figure}[ht]
            \centering
            \adjustbox{max width=0.9\textwidth}{    \begin{tikzcd}
        & [16pt] \alpha X \otimes Y 
            \arrow[r, "\id \otimes f"] 
            \arrow[d, swap, "\id \otimes \xi^\alpha_Y"]
                \arrow[dr, phantom, "\scriptstyle{\color{gray}(1)}"]
        & [16pt] \alpha X \otimes Z 
            \arrow[d, "\id \otimes \xi^\alpha_Z"]
        & [16pt] \\
        & \alpha X \otimes \alpha\alpha^{-1}Y 
            \arrow[r, "\id \otimes \alpha\alpha^{-1}(f)"] 
            \arrow[d, swap, "\eta^\alpha_{X, \alpha^{-1}Y}"]
                \arrow[dr, phantom, "\scriptstyle{\color{gray}(2)}"]
        & \alpha X \otimes \alpha\alpha^{-1}Z
            \arrow[d, "\eta^\alpha_{X, \alpha^{-1}Z}"]\\
        {}
        & \alpha (X \otimes \alpha^{-1} Y) 
            \arrow[r, "\alpha(\id \otimes \alpha^{-1}f)"] 
            \arrow[d, swap, "\alpha(\psi_{\alpha^{-1}Y})"]
                \arrow[l, phantom, ""{coordinate, name=A}]
                \arrow[l, near start, phantom, "\scriptstyle{\color{gray}(6)}"]
                \arrow[dr, phantom, "\scriptstyle{\color{gray}(3)}"]
        & \alpha (X \otimes \alpha^{-1} Z)
            \arrow[d, "\alpha(\psi_{\alpha^{-1}Z})"]
                \arrow[r, phantom, ""{coordinate, name=B}]
                \arrow[r, near start, phantom, "\scriptstyle{\color{gray}(7)}"]
        & {}\\
        & \alpha(\alpha^{-1} Y \otimes X)
            \arrow[r, "\alpha(\alpha^{-1}f \otimes \id)"] 
            \arrow[d, swap, "(\eta^\alpha_{\alpha^{-1}Y, X})^{-1}"]
                \arrow[dr, phantom, "\scriptstyle{\color{gray}(4)}"]
        & \alpha(\alpha^{-1} Z \otimes X)
            \arrow[d, "(\eta^\alpha_{\alpha^{-1}Z, X})^{-1}"]\\
        & \alpha \alpha^{-1} Y \otimes \alpha X
            \arrow[r, "\alpha\alpha^{-1}f \otimes \id"] 
            \arrow[d, swap, "(\xi^\alpha_Y)^{-1} \otimes \id"]
                \arrow[dr, phantom, "\scriptstyle{\color{gray}(5)}"]
        & \alpha\alpha^{-1} Z \otimes \alpha X
            \arrow[d, "(\xi^\alpha_Z)^{-1} \otimes \id"]\\
        & Y \otimes \alpha X 
            \arrow[r, swap, "f \otimes \id"]
            \arrow[uuuuu, <-, swap, rounded corners, to path={ -- ([xshift=-2ex]\tikztostart.west) 
                -| (A) node[pos=0.9,left]{$\scriptstyle \psi^\alpha_Y$} \tikztonodes
                |- ([xshift=-2ex]\tikztotarget.west) 
                -- (\tikztotarget)}]
        & Z \otimes \alpha X
            \arrow[uuuuu, <-, swap, rounded corners, to path={ -- ([xshift=2ex]\tikztostart.east) 
                -| (B) node[pos=0.9,right]{$\scriptstyle \psi^\alpha_Z$} \tikztonodes
                |- ([xshift=2ex]\tikztotarget.east) 
                -- (\tikztotarget)}]
    \end{tikzcd}}
            \caption{A diagram verifying the naturality of $\psi^\alpha$.}
            \label{braid-nat}
        \end{figure}
        Cells 1 and 5 commute by the naturality of $\xi^\alpha$. Cells 2 and 4 commute by the distributive naturality of the tensorator. Cell 3 commutes by the functoriality of $\alpha$ and the naturality of $\psi$. Cells 6 and 7 commute by definition of the half-braiding $\psi^\alpha$. Thus, the diagram is commutative and $\psi^\alpha$ satisfies equation \eqref{br-nat}. Similarly, we may show that $\tilde \alpha$ is well defined on morphisms in $\mc Z(\mc C)$. Consider now the diagram in figure \ref{braid-hex}:

        \begin{figure}[ht]
            \centering
            \adjustbox{max width=0.9\textwidth}{    \begin{tikzcd}[column sep=huge]
        & [12pt] \alpha X \otimes Y \otimes Z 
            \arrow[dd, swap, "\id \otimes \xi^\alpha_Y \otimes \id"]
            \arrow[dr, "\id \otimes \xi^\alpha_Z"]
            \arrow[ddrr, bend left=20, "\id \otimes \xi^\alpha_{Y \otimes Z}"]
                \arrow[ddr, bend right=5, phantom, "\scriptstyle{\color{gray}(1)}"]
        \\
        && \alpha X \otimes Y \otimes \alpha\alpha^{-1} Z 
            \arrow[d, swap, "\id \otimes \xi^\alpha_Y \otimes \id"]
                \arrow[dr, phantom, near start, "\scriptstyle{\color{gray}(2)}"]
        \\
        & \alpha X \otimes \alpha\alpha^{-1}Y \otimes Z
            \arrow[d, swap, "\eta^\alpha_{X, \alpha^{-1}Y} \otimes \id"]
            \arrow[r, "\id \otimes \xi^\alpha_Z"]
                \arrow[dr, phantom, "\scriptstyle{\color{gray}(3)}"]
        & \alpha X \otimes \alpha\alpha^{-1}Y \otimes \alpha\alpha^{-1} Z
            \arrow[d, swap, "\eta^\alpha_{X, \alpha^{-1}Y} \otimes \id"]
            \arrow[dr, bend left=5, swap, near start, "\id \otimes \eta^\alpha_{\alpha^{-1}Y, \alpha^{-1}Z}"]
            \arrow[r, "\id \otimes \eta^{\alpha\alpha^{-1}}_{Y, Z}"]
                \arrow[dr, phantom, near end, bend left=15, "\scriptstyle{\color{gray}(5)}"]
        & \alpha X \otimes \alpha\alpha^{-1}(Y \otimes Z)
            \arrow[d, "\id \otimes \alpha(\eta^{\alpha^{-1}}_{Y, Z})^{-1}"]
            \arrow[ddr, bend left, "\eta^\alpha_{X, \alpha^{-1}(Y \otimes Z)}"]
                \arrow[ddr, phantom, "\scriptstyle{\color{gray}(6)}"]
            \\
        {}
        & \alpha (X \otimes \alpha^{-1} Y) \otimes Z
            \arrow[d, swap, "\alpha(\psi_{\alpha^{-1}Y}) \otimes \id"]
            \arrow[r, "\id \otimes \xi^\alpha_Z"]
                \arrow[l, phantom, ""{coordinate, name=A}]
                \arrow[dr, phantom, "\scriptstyle{\color{gray}(7)}"]
                \arrow[l, phantom, near start, "\scriptstyle{\color{gray}(23)}"]
        & \alpha (X \otimes \alpha^{-1} Y) \otimes \alpha\alpha^{-1} Z
            \arrow[d, swap, "\alpha(\psi_{\alpha^{-1}Y}) \otimes \id"]
            \arrow[dr, bend left=5, swap, near start, "\eta^\alpha_{X \otimes \alpha^{-1}Y, \alpha^{-1}Z}"]
                \arrow[r, phantom, "\scriptstyle{\color{gray}(4)}"]
        & \alpha X \otimes \alpha(\alpha^{-1} Y \otimes \alpha^{-1} Z)
            \arrow[d, "\eta^\alpha_{X, \alpha^{-1}Y \otimes \alpha^{-1}Z}"]
        \\
        & \alpha(\alpha^{-1} Y \otimes X) \otimes Z
            \arrow[d, swap, "(\eta^\alpha_{\alpha^{-1}Y, X})^{-1} \otimes \id"]
            \arrow[r, "\id \otimes \xi^\alpha_Z"]
                \arrow[ddr, bend left=5, near start, phantom, "\scriptstyle{\color{gray}(9)}"]
        & \alpha(\alpha^{-1} Y \otimes X) \otimes \alpha\alpha^{-1} Z
            \arrow[dd, swap, near start, "(\eta^\alpha_{\alpha^{-1}Y, X})^{-1} \otimes \id"]
            \arrow[ddr, bend left=15, near start, swap, "\eta^\alpha_{\alpha^{-1}Y \otimes X, \alpha^{-1}Z}"]
                \arrow[r, phantom, "\scriptstyle{\color{gray}(8)}"]
        & \alpha(X \otimes \alpha^{-1} Y \otimes \alpha^{-1} Z)
            \arrow[dd, near end, "\alpha(\psi_{\alpha^{-1}Y} \otimes \id)"]
            \arrow[r, swap, "\alpha(\id \otimes \eta^{\alpha^{-1}}_{Y, Z})"]
        & \alpha(X \otimes \alpha^{-1}(Y \otimes Z))
            \arrow[dddd, bend left, near end, "\alpha(\psi_{\alpha^{-1}(Y \otimes Z)})"]
            \\
        & \alpha \alpha^{-1} Y \otimes \alpha X \otimes Z
            \arrow[d, swap, "(\xi^\alpha_Y)^{-1} \otimes \id"]
            \arrow[dr, "\id \otimes \xi^\alpha_Z"]
            \\
        & Y \otimes \alpha X \otimes Z
            \arrow[d, swap, "\id \otimes \xi^\alpha_Z"]
            \arrow[uuuuuu, <-, swap, rounded corners, to path={ -- ([xshift=-2ex]\tikztostart.north west) 
                -| (A) node[pos=0.9,left]{$\scriptstyle \psi^\alpha_Y$} \tikztonodes
                |- ([xshift=-2ex]\tikztotarget.west) 
                -- (\tikztotarget)}] 
                \arrow[r, phantom, "\scriptstyle{\color{gray}(10)}"]
        & \alpha \alpha^{-1} Y \otimes \alpha X \otimes \alpha\alpha^{-1}Z
            \arrow[dl, "(\xi^\alpha_Y)^{-1} \otimes \id"]
            \arrow[dd, swap, near end, "\id \otimes \eta^\alpha_{X, \alpha^{-1}Z}"]
                \arrow[r, phantom, "\scriptstyle{\color{gray}(12)}"]
        & \alpha(\alpha^{-1} Y \otimes X \otimes \alpha^{-1} Z)
            \arrow[ddl, bend left=15, swap, near end, "(\eta^\alpha_{\alpha^{-1}Y, X \otimes \alpha^{-1}Z})^{-1}"]
            \arrow[dd, near start, "\alpha(\id \otimes \psi_{\alpha^{-1}Z})"]
                \arrow[r, phantom, ""{coordinate, name=C}]
                \arrow[rr, phantom, ""{coordinate, name=D}]
                \arrow[r, phantom, near start, "\scriptstyle{\color{gray}(13)}"]
                \arrow[rr, phantom, "\scriptstyle{\color{gray}(14)}"]
                \arrow[rr, phantom, near end, "\scriptstyle{\color{gray}(25)}"]
        & {} & {}
            \\
        & Y \otimes \alpha X \otimes \alpha\alpha^{-1}Z 
            \arrow[d, swap, "\id \otimes \eta^\alpha_{X, \alpha^{-1}Z}"]
            \\
        {}
        & Y \otimes \alpha (X \otimes \alpha^{-1} Z) 
            \arrow[d, swap, "\id \otimes \alpha(\psi_{\alpha^{-1}Z})"]
            \arrow[<-, r, "(\xi^\alpha_Y)^{-1} \otimes \id"]
                \arrow[l, phantom, ""{coordinate, name=B}]
                \arrow[uur, bend left=5, near start, phantom, "\scriptstyle{\color{gray}(11)}"]
                \arrow[dr, phantom, "\scriptstyle{\color{gray}(15)}"]
                \arrow[l, phantom, near start, "\scriptstyle{\color{gray}(24)}"]
        & \alpha\alpha^{-1} Y \otimes \alpha(X \otimes \alpha^{-1} Z)
            \arrow[d, swap, "\id \otimes \alpha(\psi_{\alpha^{-1}Z})"]
                \arrow[r, phantom, "\scriptstyle{\color{gray}(16)}"]
        & \alpha(\alpha^{-1} Y \otimes \alpha^{-1} Z \otimes X)
            \arrow[d, "(\eta^\alpha_{\alpha^{-1}Y \otimes \alpha^{-1}Z, X})^{-1}"]
            \arrow[r, "\alpha(\eta^{\alpha^{-1}}_{Y, Z} \otimes \id)"]
            \arrow[uuuu, <-, swap, rounded corners, to path={ -- ([xshift=-2ex, yshift=4ex]\tikztostart.north east) 
                -| ([xshift=-2ex]C) node[pos=0.85,right]{$\scriptstyle \alpha(\psi_{\alpha^{-1}Y \otimes \alpha^{-1}Z})$} \tikztonodes
                |- ([xshift=-2ex, yshift=-4ex]\tikztotarget.south east) 
                -- (\tikztotarget)}] 
        & \alpha(\alpha^{-1}(Y \otimes Z) \otimes X)
            \arrow[ddl, bend left, "(\eta^\alpha_{\alpha^{-1}(Y \otimes Z), X})^{-1}"]
        \\
        & Y \otimes \alpha(\alpha^{-1} Z \otimes X)
            \arrow[d, swap, "\id \otimes (\eta^\alpha_{\alpha^{-1}Z, X})^{-1}"]
            \arrow[<-, r, "(\xi^\alpha_Y)^{-1} \otimes \id"]
                \arrow[dr, phantom, "\scriptstyle{\color{gray}(17)}"]
        & \alpha\alpha^{-1} Y \otimes \alpha(\alpha^{-1} Z \otimes X)
            \arrow[d, swap, "\id \otimes (\eta^\alpha_{\alpha^{-1}Z, X})^{-1}"]
            \arrow[<-, ur, bend right=5, "(\eta^\alpha_{\alpha^{-1}Y, X \otimes \alpha^{-1}Z})^{-1}"]
                \arrow[r, phantom, "\scriptstyle{\color{gray}(18)}"]
        & \alpha(\alpha^{-1} Y \otimes \alpha^{-1} Z) \otimes \alpha X
            \arrow[d, "\alpha(\eta^{\alpha^{-1}}_{Y, Z}) \otimes \id"]
            \\
        & Y \otimes \alpha \alpha^{-1} Z \otimes \alpha X
            \arrow[dd, swap, "\id \otimes (\xi^\alpha_Z)^{-1} \otimes \id"]
            \arrow[<-, r, "(\xi^\alpha_Y)^{-1} \otimes \id"]
        & \alpha \alpha^{-1} Y \otimes \alpha \alpha^{-1} Z \otimes \alpha X
            \arrow[d, swap, "\id \otimes (\xi^\alpha_Z)^{-1} \otimes \id"]
            \arrow[<-, r, swap, "\eta^{\alpha\alpha^{-1}}_{Y, Z} \otimes \id"]
            \arrow[<-, ur, bend right=5, "(\eta^\alpha_{\alpha^{-1}Y, \alpha^{-1}Z})^{-1} \otimes \id"]
                \arrow[ur, phantom, near end, bend right=15, "\scriptstyle{\color{gray}(19)}"]
        & \alpha\alpha^{-1} (Y \otimes Z) \otimes X
            \arrow[ddll, bend left=20, "(\xi^\alpha_{Y \otimes Z})^{-1} \otimes \id"]
                \arrow[uur, phantom, "\scriptstyle{\color{gray}(20)}"]
            \\
        && \alpha \alpha^{-1} Y \otimes Z \otimes \alpha X
            \arrow[dl, "(\xi^\alpha_Y)^{-1} \otimes \id"]
                \arrow[ur, phantom, near start, "\scriptstyle{\color{gray}(22)}"]
            \\
        & Y \otimes Z \otimes \alpha X 
            \arrow[uuuuuu, <-, swap, rounded corners, to path={ -- ([xshift=-2ex]\tikztostart.west) 
                -| (B) node[pos=0.9,left]{$\scriptstyle \psi^\alpha_Y$} \tikztonodes
                |- ([xshift=-2ex]\tikztotarget.south west) 
                -- (\tikztotarget)}] 
            \arrow[uuuuuuuuuuuu, <-, swap, rounded corners, to path={ -- ([xshift=-2ex, yshift=-2ex]\tikztostart.south east) 
                -| ([xshift=20ex]D) node[pos=1.1,right]{$\scriptstyle \psi^\alpha_{Y\otimes Z}$} \tikztonodes
                |- ([xshift=-2ex, yshift=2ex]\tikztotarget.north east) 
                -- (\tikztotarget)}] 
                    \arrow[uur, bend left=5, phantom, "\scriptstyle{\color{gray}(21)}"]
    \end{tikzcd}}
            \caption{A diagram verifying $\psi^\alpha$ satisfies the hexagon equations.}
            \label{braid-hex}
        \end{figure}
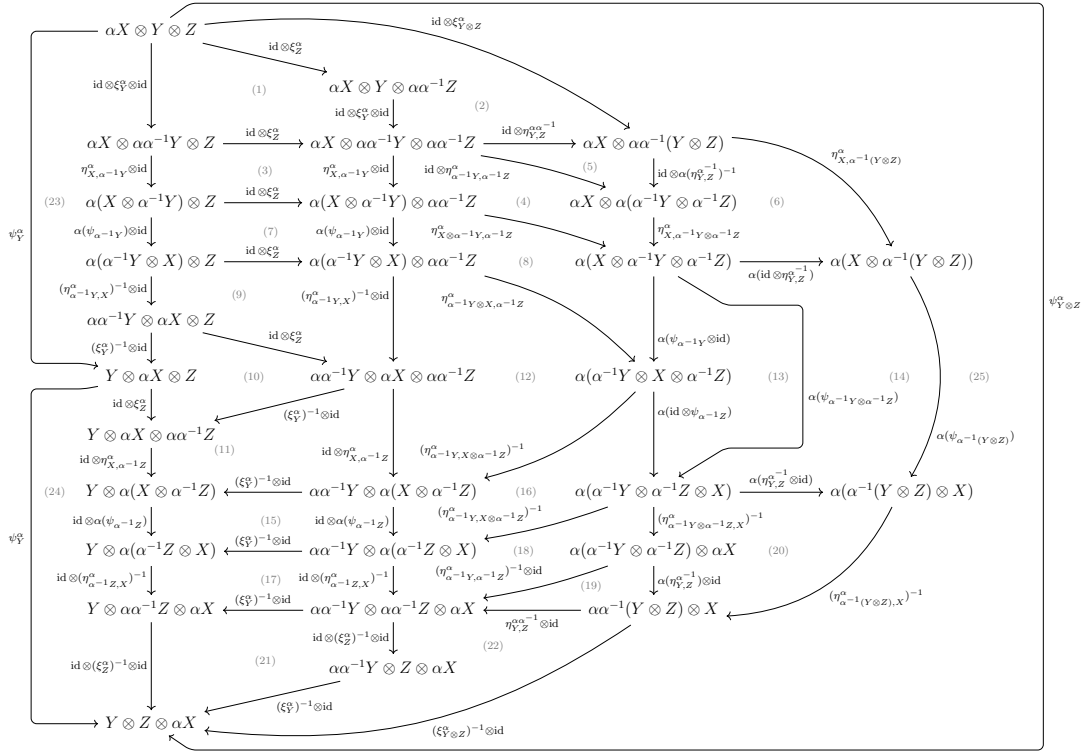

        Cells 1, 3, 7, 9, 10, 11, 15, 17, and 21 commute by the naturality of the monoidal product. Cells 2 and 22 commute as $\xi$ is a monoidal natural isomorphism. Cells 4, 12 and 18 commute by the associative naturality of the tensorator. Cells 5 and 19 commute by the definition of the composition of monoidal functors. Cells 6, 8, 16 and 20 commute by the distributive naturality of the tensorator. Cell 13 commutes by the functoriality of $\alpha$ and since $\psi$ satisfies \eqref{br-hex}. Cell 14 commutes by the functoriality of $\alpha$ and since $\psi$ satisfies \eqref{br-nat}. Cells 23, 24 and 25 commute by the definition of $\psi^\alpha$. Thus, the diagram commutes and $\psi^\alpha$ satisfies equation \eqref{br-hex}. Therefore, $\psi^\alpha$ is a half-braiding and $\tilde\alpha$ is well-defined on objects of $\mc Z(\mc C)$.

        We now demonstrate that $(\tilde\alpha, \tilde\eta^\alpha)$ is a monoidal functor. As $(\alpha, \eta^\alpha)$ is a monoidal functor, we need only show that the tensorator is compatible with the monoidal product in the center. Let $(X, \psi), (Y, \varphi) \in \mc Z(\mc C)$ and denote the braiding of their product $\bar \psi := (\psi \otimes \id) \circ (\id \otimes \varphi)$. We thus need that that:
        \begin{align} \label{monoidal}
            (\id_Z \otimes \, \eta^\alpha_{X, Y})
            \circ (\psi^\alpha_Z \otimes \id_{\alpha Y}) 
            \circ (\id_{\alpha X} \otimes \, \varphi^\alpha_Z) 
            = \bar \psi^\alpha_Z \circ (\eta^{\alpha}_{X, Y} \otimes \id_Z), \qquad \text{ for all } Z \in \mc C.
        \end{align}
        
        Consider diagram \ref{aut-mon}.
        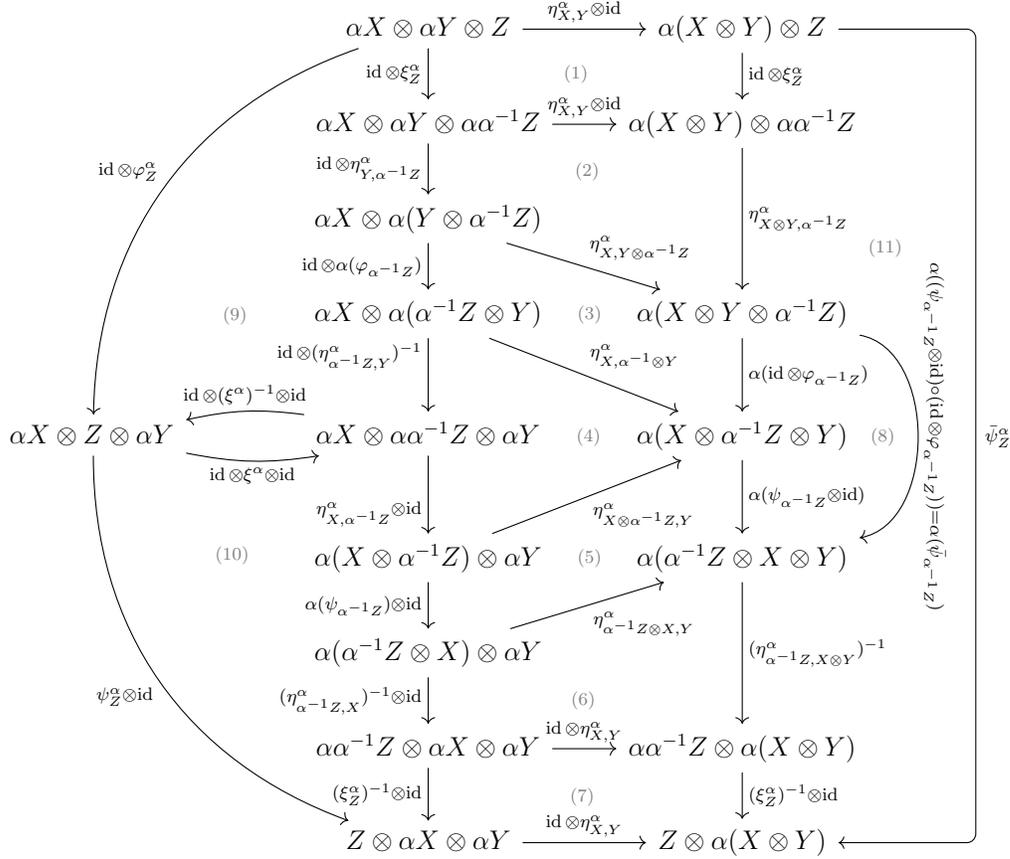
\begin{figure}[ht]
            \centering
            \adjustbox{max width=0.9\textwidth}{    \begin{tikzcd}
        & [24pt] \alpha X \otimes \alpha Y \otimes Z 
            \arrow[r, "\eta^\alpha_{X, Y} \otimes \id"]
            \arrow[d, swap, "\id \otimes \xi^{\alpha}_Z"]
            \arrow[ddddl, bend right=35, swap, "\id \otimes \varphi^\alpha_Z"]
                \arrow[dr, phantom, "\scriptstyle{\color{gray}(1)}"]
        & \alpha(X \otimes Y) \otimes Z 
            \arrow[d, "\id \otimes \xi^{\alpha}_Z"]
        & [72pt] \\
        & \alpha X \otimes \alpha Y \otimes \alpha\alpha^{-1} Z 
            \arrow[r, "\eta^\alpha_{X, Y} \otimes \id"]
            \arrow[d, swap, "\id \otimes \eta^\alpha_{Y, \alpha^{-1}Z}"]
        & \alpha(X \otimes Y) \otimes \alpha\alpha^{-1}Z 
            \arrow[dd, "\eta^\alpha_{X \otimes Y, \alpha^{-1}Z}"]
                \arrow[dddr, pos=0.35, phantom, "\scriptstyle{\color{gray}(11)}"]\\
        & \alpha X \otimes \alpha(Y \otimes \alpha^{-1} Z) 
            \arrow[d, swap, "\id \otimes \alpha (\varphi_{\alpha^{-1}Z})"]
            \arrow[dr, "\eta^\alpha_{X, Y \otimes \alpha^{-1}Z}"]
                \arrow[ur, phantom, "\scriptstyle{\color{gray}(2)}"]
                \arrow[ddl, bend right=10, phantom, "\scriptstyle{\color{gray}(9)}"]\\
        & \alpha X \otimes \alpha(\alpha^{-1}Z \otimes Y) 
            \arrow[d, swap, near start, "\id \otimes (\eta^\alpha_{\alpha^{-1}Z, Y})^{-1}"] 
            \arrow[dr, "\eta^\alpha_{X, \alpha^{-1} \otimes Y}"]
                \arrow[r, phantom, "\scriptstyle{\color{gray}(3)}"]
        & \alpha(X \otimes Y \otimes \alpha^{-1}Z) 
            \arrow[d, "\alpha(\id \otimes \varphi_{\alpha^{-1}Z})"]
            \arrow[dd, bend left=80, "\alpha((\psi_{\alpha^{-1}Z} \otimes \id)\circ(\id \otimes \varphi_{\alpha^{-1}Z})) = \alpha(\bar \psi_{\alpha^{-1}Z})" {anchor=south, rotate=-90}] 
            \\ [12pt]
        \alpha X \otimes Z \otimes \alpha Y 
            \arrow[r, bend right=10, swap, "\id \otimes \xi^{\alpha} \otimes \id"] 
            \arrow[ddddr, bend right=35, swap, "\psi^\alpha_Z \otimes \id"]
        & \alpha X \otimes \alpha \alpha^{-1} Z \otimes \alpha Y 
            \arrow[d, swap, near end, "\eta^\alpha_{X, \alpha^{-1}Z} \otimes \id"] 
            \arrow[l, swap, bend right=10, "\id \otimes (\xi^{\alpha})^{-1} \otimes \id"]
                \arrow[r, phantom, "\scriptstyle{\color{gray}(4)}"]
        & \alpha(X \otimes \alpha^{-1}Z \otimes Y) 
            \arrow[d, "\alpha(\psi_{\alpha^{-1}Z} \otimes \id)"] 
                \arrow[r, phantom, ""{coordinate, name=Z}] 
                \arrow[r, phantom, pos=0.1, "\scriptstyle{\color{gray}(8)}"]
        & {}\\ [12pt]
        & \alpha(X \otimes \alpha^{-1}Z) \otimes \alpha Y 
            \arrow[d, swap, "\alpha (\psi_{\alpha^{-1}Z}) \otimes \id"]
            \arrow[ur, swap, "\eta^\alpha_{X \otimes \alpha^{-1}Z, Y}"]
                \arrow[r, phantom, "\scriptstyle{\color{gray}(5)}"]
        & \alpha(\alpha^{-1}Z \otimes X \otimes Y)
            \arrow[dd, "(\eta^\alpha_{\alpha^{-1}Z, X \otimes Y})^{-1}"]\\
        & \alpha(\alpha^{-1}Z \otimes X) \otimes \alpha Y
            \arrow[d, swap, "(\eta^\alpha_{\alpha^{-1}Z, X})^{-1} \otimes \id"] 
            \arrow[ur, swap, "\eta^\alpha_{\alpha^{-1}Z \otimes X, Y}"]
                \arrow[dr, phantom, "\scriptstyle{\color{gray}(6)}"]
                \arrow[uul, bend left=10, phantom, "\scriptstyle{\color{gray}(10)}"]\\
        & \alpha \alpha^{-1}Z \otimes \alpha X \otimes \alpha Y 
            \arrow[r, "\id \otimes \eta^\alpha_{X, Y}"]
            \arrow[d, swap, "(\xi^\alpha_Z)^{-1} \otimes \id"]
        & \alpha\alpha^{-1}Z\otimes\alpha(X \otimes Y) 
            \arrow[d, "(\xi^{\alpha}_Z)^{-1} \otimes \id"]\\
        & Z \otimes \alpha X \otimes \alpha Y 
            \arrow[r, "\id \otimes \eta^\alpha_{X, Y}"]
                \arrow[ur, phantom, "\scriptstyle{\color{gray}(7)}"]
        & Z \otimes \alpha(X \otimes Y)
            \arrow[uuuuuuuu, <-, swap, rounded corners, to path={ -- ([xshift=2ex]\tikztostart.east) 
            -| (Z) node[pos=1,right]{$\scriptstyle \bar \psi^\alpha_Z$} \tikztonodes
            |- ([xshift=2ex]\tikztotarget.east) 
            -- (\tikztotarget)}]
    \end{tikzcd}}.
            \caption{A diagram verifying that $(\tilde\alpha, \tilde\eta^\alpha)$ is monoidal.}
            \label{aut-mon}
        \end{figure}
        Cell 1 and 7 commute by the naturality of the monoidal product. Cells 2, 4, and 6 commute by the associative naturality of the tensorator. Cells 3 and 5 commute by the distributive naturality of the tensorator. Cell 8 commutes by the functoriality of $\alpha$. Cells 9, 10 and 11 commute by the definition of the half-braidings.
        
        We now verify that $\underline{I^L}$ is well-defined on 2-morphisms. Given $\pi \in \Hom((\alpha, \eta^\alpha), (\beta, \eta^\alpha))$ a monoidal natural transformation, we define $\underline{I^L}(\pi)_X$ to be the lift of $\pi_{\Forg X}$. By abuse of notation, we write $\pi_{\Forg X} = \pi_X$. For this lift to be well-defined, we need to show that $\pi_X$ is a morphism in the center. That is, given $(X, \psi) \in \mc Z(\mc C)$, for any $Y \in \mc C$:
        \begin{align*} 
            \psi^\beta_Y \circ (\pi_X \otimes \id_Y) = (\id_Y \otimes \pi_X) \circ \psi^\alpha_Y.
        \end{align*} 

        Consider diagram \ref{nat-cent}.
        \begin{figure}[ht]
            \centering
            \adjustbox{max width=0.9\textwidth}{    \begin{tikzcd} 
        &&[12pt]&[12pt] \alpha X \otimes Y
            \arrow[d, swap, "\id \otimes \xi^\alpha_Y"]
            \arrow[ddddr, bend left, "\psi^\alpha_Y"] \\
        \alpha X \otimes Y 
            \arrow[r, "\id \otimes \xi^\beta_Y"]
            \arrow[d, swap, "\pi_X \otimes \id"]
            \arrow[urrr, bend left=20, "\id \otimes \left((\xi^\alpha_Y)^{-1} \circ\, \alpha(\pi^{-1}_Y) \,\circ\, \pi^{-1}_{\beta^{-1}Y}\circ\, \xi^\beta_Y\right)"]
                \arrow[dr, phantom, "\scriptstyle{\color{gray}(1)}"]
        & \alpha X \otimes \beta \beta^{-1} Y
            \arrow[d, "\pi_X \otimes \id"]
            \arrow[r, "\id \otimes \pi^{-1}_{\beta^{-1}Y}"]
                \arrow[dr, phantom, "\scriptstyle{\color{gray}(2)}"]
                \arrow[urr, bend left=5, phantom, "\scriptstyle{\color{gray}(10)}"]
        & \alpha X \otimes \alpha \beta^{-1} Y
            \arrow[d, "\eta^\alpha_{X, \beta^{-1}Y}"] 
            \arrow[r, "\id \otimes \alpha(\pi^{-1}_Y)"]
                \arrow[dr, phantom, "\scriptstyle{\color{gray}(3)}"]
        & \alpha X \otimes \alpha \alpha^{-1} Y
            \arrow[d, "\eta^\alpha_{X, \alpha^{-1}Y}"]\\
        \beta X \otimes Y 
            \arrow[r, swap, "\id \otimes \xi^\beta_Y"] 
            \arrow[ddddr, bend right, swap, "\psi^\beta_Y"]
        & \beta X \otimes \beta \beta^{-1} Y
            \arrow[d, "\eta^\beta_{X, \beta^{-1}Y}"]
        & \alpha(X \otimes \beta^{-1} Y)
            \arrow[d, "\alpha(\psi_{\beta^{-1}Y})"] 
            \arrow[r, "\alpha(\id \otimes \pi^{-1}_Y)"]
                \arrow[dr, phantom, "\scriptstyle{\color{gray}(5)}"]
        & \alpha(X \otimes \alpha^{-1} Y)
            \arrow[d, "\alpha(\psi_{\alpha^{-1}Y})"]\\
        & \beta(X\otimes \beta^{-1}Y)
            \arrow[d, "\beta(\psi_{\beta^{-1}Y})"]
            \arrow[ur, swap, "\pi_{X \otimes \beta^{-1}Y}^{-1}"]
                \arrow[r, phantom, "\scriptstyle{\color{gray}(4)}"]
        & \alpha(\beta^{-1}Y \otimes X)
            \arrow[d, "(\eta^\alpha_{\beta^{-1}Y, X})^{-1}"] 
            \arrow[r, "\alpha(\pi^{-1}_Y \otimes \id)"]
                \arrow[dr, bend left=5, phantom, "\scriptstyle{\color{gray}(6)}"]
        & \alpha(\alpha^{-1}Y \otimes X)
            \arrow[d, "(\eta^{\alpha}_{\alpha^{-1}Y, X})^{-1}"]\\
        & \beta(\beta^{-1}Y \otimes X)
            \arrow[d, "(\eta^\beta_{\beta^{-1}Y, X)})^{-1}"]
            \arrow[ur, swap, "\pi_{\beta^{-1}Y \otimes X}^{-1}"]
                \arrow[uul, bend left=5,  phantom, "\scriptstyle{\color{gray}(12)}"]
        & \alpha\beta^{-1}Y \otimes \alpha X
            \arrow[d, "\id \otimes \pi_X"]
            \arrow[r, "\alpha(\pi^{-1}_Y) \otimes \id"]
                \arrow[dr, phantom, "\scriptstyle{\color{gray}(8)}"]
        & \alpha \alpha^{-1}Y \otimes \alpha X
            \arrow[d, "\id \otimes \pi_X"]
            \arrow[r, "(\xi^\alpha_Y)^{-1} \otimes \id"]
                \arrow[dr, phantom, "\scriptstyle{\color{gray}(9)}"]
        & Y \otimes \alpha X
            \arrow[d, "\id \otimes \pi_X"]
                \arrow[uul, bend right=5, phantom, "\scriptstyle{\color{gray}(13)}"]\\
        & \beta\beta^{-1}Y \otimes \beta X
            \arrow[d, "(\xi^\beta_Y)^{-1}\otimes \id"] 
            \arrow[r, "(\pi_{\beta^{-1}Y})^{-1}\otimes \id"]
                \arrow[ur, near end, phantom, "\scriptstyle{\color{gray}(7)}"]
        & \alpha\beta^{-1}Y \otimes \beta X
            \arrow[r, "\alpha(\pi^{-1}_Y) \otimes \id"]
        & \alpha \alpha^{-1} Y \otimes \beta X
            \arrow[r, "(\xi^\alpha_Y)^{-1} \otimes \id"]
        & Y \otimes \beta X
            \arrow[dlll, bend left=20, "\left((\xi^\alpha_Y)^{-1} \circ \, \alpha(\pi^{-1}_Y) \,\circ \,\pi^{-1}_{\beta^{-1}Y}\circ\, \xi^\beta_Y\right)^{-1} \otimes \id"]\\
        & Y \otimes \beta X
            \arrow[urr, bend right=5, phantom, "\scriptstyle{\color{gray}(11)}"]
    \end{tikzcd}}.
            \caption{A diagram verifying the lift of a natural transformation is central.}
            \label{nat-cent}
        \end{figure}
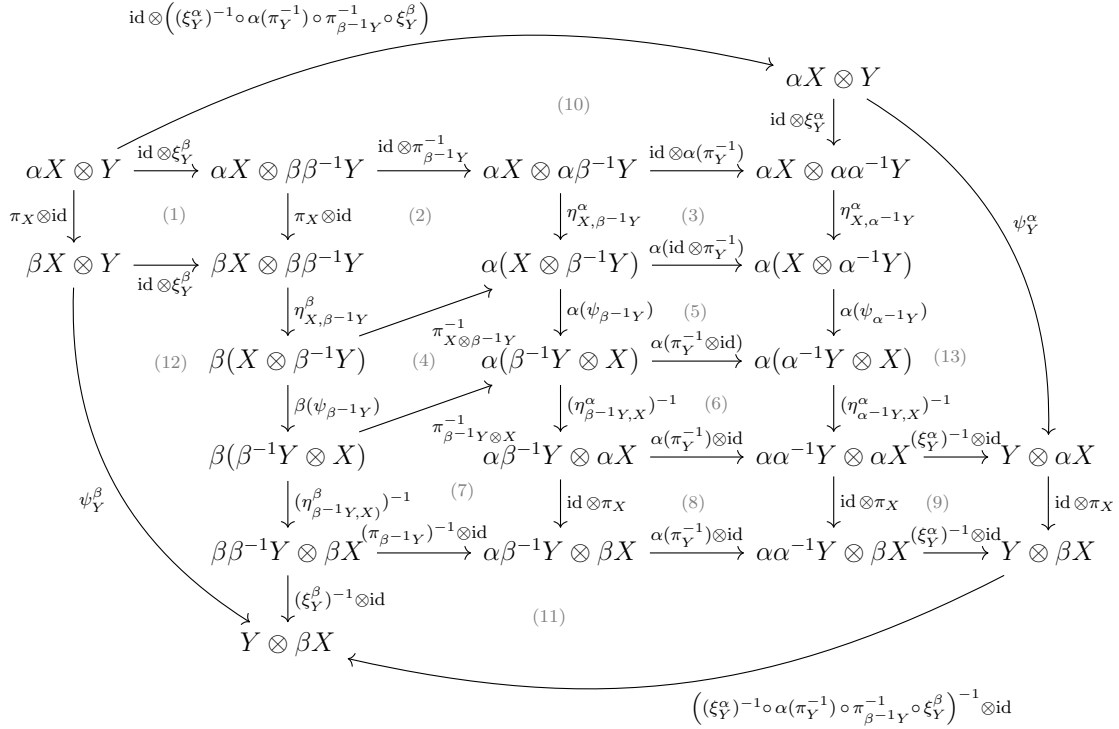
        Cells 1, 8 and 9 commute by the naturality of the monoidal product. Cells 2 and 7 commute as $\pi$ is a monoidal natural isomorphism. Cells 3 and 6 commute by the distributive naturality of the tensorator. Cells 4 and 5 commute by the naturality of $\pi$. Cells 12 and 13 commute by the definition of the half-braiding. Cells 10 and 11 commute by definition.
        
        Thus, the diagram is commutative and we obtain:
        \begin{align*}
            & \psi_Y^\beta \circ (\pi_X \otimes \id) \\
            & = (((\xi^\alpha_Y)^{-1} \circ \alpha(\pi^{-1}_Y) \circ \pi^{-1}_{\beta^{-1}Y}\circ \xi^\beta_Y)^{-1} \otimes \id) \\
            & \quad \circ (\id \otimes\, \pi_X) \circ \psi_Y^\alpha \circ (\id \otimes ((\xi^\alpha_Y)^{-1} \circ \alpha(\pi^{-1}_Y) \circ \pi^{-1}_{\beta^{-1}Y}\circ \xi^\beta_Y)) \\
            & = ((((\xi^\alpha_Y)^{-1} \circ \alpha(\pi^{-1}_Y) \circ \pi^{-1}_{\beta^{-1}Y}\circ \xi^\beta_Y)^{-1} \\
            & \quad \circ ((\xi^\alpha_Y)^{-1} \circ \alpha(\pi^{-1}_Y) \circ \pi^{-1}_{\beta^{-1}Y}\circ \xi^\beta_Y)) \otimes \id) \circ (\id \otimes\, \pi_X) \circ \psi_Y^\alpha \\
            & = (\id \otimes \id) \circ (\id \otimes\, \pi_X) \circ \psi^\alpha_Y = (\id \otimes\, \pi_X) \circ \psi^\alpha_Y,
        \end{align*}
        as $\psi_Y^\alpha$ is a half-braiding. We have thus shown that $\pi_X$ is a morphism in the center. 

        Recall that $\mc Z(\mc C)$ is a braided category with the braiding $C_{(X, \psi), (Y, \phi)} = \psi_Y$. We now show that $(\tilde\alpha, \eta^\alpha)$ is \emph{braided} monoidal endofunctor. For $\tilde\alpha$ to be braided, we need that:
        \begin{align*}
            \alpha(C_{(X, \psi), (Y, \phi)}) \circ \eta^\alpha_{X, Y} = \eta^\alpha_{Y, X} \circ C_{\tilde\alpha(X, \psi), \tilde\alpha(Y, \phi)}.
        \end{align*}
        Consider the morphism $\alpha(\xi^{\alpha^{-1}}_Y)^{-1} \circ \xi^\alpha_{\alpha Y}$, and their composition are morphisms $Y \xrightarrow{} Y$. As $C_{\tilde\alpha(X, \psi), \tilde\alpha(Y, \phi)} = \psi^\alpha_{\alpha Y}$ is a half-braiding, we have:
        \begin{align*}
            & C_{\tilde\alpha(X, \psi), \tilde\alpha(Y, \phi)} = \psi_{\alpha Y}^\alpha = \psi_{\alpha Y}^\alpha \circ (\id \otimes \id) \\
            & = \psi_{\alpha Y}^\alpha \circ (\id \otimes (\alpha(\xi^{\alpha^{-1}}_Y)^{-1} \circ \xi^\alpha_{\alpha Y} \circ (\xi^\alpha_{\alpha Y})^{-1} \circ \alpha(\xi^{\alpha^{-1}}_Y)))) \\
            & = ((\alpha(\xi^{\alpha^{-1}}_Y)^{-1} \circ \xi^\alpha_{\alpha Y}) \otimes \id) \circ \psi_{\alpha Y}^\alpha \circ (\id \otimes ((\xi^\alpha_{\alpha Y})^{-1} \circ \alpha(\xi^{\alpha^{-1}}_Y))) \\
            & = ((\alpha(\xi^{\alpha^{-1}}_Y)^{-1} \circ \xi^\alpha_{\alpha Y}) \otimes \id) \, \circ \\
            & \qquad ((\xi^\alpha_{\alpha Y})^{-1} \otimes \id) \circ (\eta^\alpha_{\alpha^{-1}\alpha Y, X})^{-1} \circ \alpha(\psi_{\alpha\alpha^{-1}Y}) \circ \eta^\alpha_{X, \alpha^{-1}\alpha Y} \circ (\id \otimes \xi^\alpha_{\alpha Y})  \,\circ \\
            & \qquad\qquad (\id \otimes ((\xi^\alpha_{\alpha Y})^{-1} \circ \alpha(\xi^{\alpha^{-1}}_Y))) \\
            & = (\alpha(\xi^{\alpha^{-1}}_Y)^{-1} \otimes \id) \circ (\eta^\alpha_{\alpha^{-1}\alpha Y, X})^{-1} \circ \alpha(\psi_{\alpha \alpha^{-1}Y}) \circ \eta^\alpha_{X, \alpha^{-1}\alpha Y} \circ (\id \otimes \alpha(\xi^{\alpha^{-1}}_Y)).
        \end{align*}
        However, by the distributive naturality of the tensorator, the functoriality of $\alpha$, and the naturality of the half-braiding $\psi$, the $\alpha(\xi^{\alpha^{-1}}_Y)$ terms pass over the braiding and we have that:
        \begin{align*}
            C_{\tilde\alpha(X, \psi), \tilde\alpha(Y, \phi)} = (\eta^\alpha_{Y, X})^{-1} \circ \alpha(\psi_Y) \circ \eta^\alpha_{X, Y}
            = (\eta^\alpha_{Y, X})^{-1} \circ \alpha(C_{(X, \psi), (Y, \phi)}) \circ \eta^\alpha_{X, Y},
        \end{align*}
        as desired. Thus, $(\tilde \alpha, \eta^\alpha)$ is a braided monoidal endofunctor.

        We now show that $(\tilde\alpha, \tilde\eta^\alpha)$ is an \emph{autoequivalence}. As $(\alpha, \eta^\alpha)$ is an autoequivalence, we need only show that $(\psi^\beta)^\alpha = \psi^{\alpha\beta}$. Consider diagram \ref{half-comp}:
        \begin{figure}[ht]
            \centering
            \adjustbox{max width=0.9\textwidth}{    \begin{tikzcd}
        & [24pt] & [12pt] \alpha\beta X \otimes Y
            \arrow[d, "\id \otimes \xi^{\alpha_Y}"]\\
        &\alpha\beta X \otimes Y
            \arrow[d, swap, "\id \otimes \xi^{\alpha\beta}_Y"] 
            \arrow[ur, "((\xi^\alpha_Y)^{-1} \circ \,\alpha(\xi^\beta_{\alpha^{-1}Y})^{-1} \circ \, \xi^{\alpha\beta}_Y) \otimes \id"]
                \arrow[r, phantom, "\scriptstyle{\color{gray}(1)}"]
        & \alpha\beta X \otimes \alpha\alpha^{-1}Y 
            \arrow[dl, swap, <-, "\id \otimes \alpha(\xi^\beta_{\alpha^{-1}Y})^{-1}"]
            \arrow[d, "\eta^\alpha_{\beta X, \alpha^{-1}Y}"]
        &[12pt] & [64pt]\\
        &\alpha \beta X \otimes \alpha \beta(\alpha \beta)^{-1} Y
            \arrow[d, swap, "\eta^{\alpha\beta}_{X, (\alpha\beta)^{-1}Y}"] 
            \arrow[dr, "\eta^\alpha_{\beta X, \beta (\alpha \beta)^{-1}Y}"]
                \arrow[r, phantom, "\scriptstyle{\color{gray}(2)}"]
                \arrow[dr, bend right=10, phantom, "\scriptstyle{\color{gray}(3)}"]
        & \alpha(\beta X \otimes \alpha^{-1} Y) 
            \arrow[d, "\alpha(\id \otimes \xi^{\beta}_{\alpha^{-1}Y})"] \\
        { }&\alpha \beta (X \otimes (\alpha \beta)^{-1} Y)
            \arrow[d, swap, "\alpha\beta(\psi_{(\alpha\beta)^{-1}Y})"]
                \arrow[l, phantom, ""{coordinate, name=C}] 
                \arrow[l, phantom, near start, "\scriptstyle{\color{gray}(8)}"]
        & \alpha(\beta X \otimes \beta \beta^{-1} \alpha^{-1} Y)
            \arrow[l, "\alpha(\eta^\beta_{X, (\alpha\beta)^{-1}Y})"] 
                \arrow[r, phantom, ""{coordinate, name=A}] 
                \arrow[rr, phantom, ""{coordinate, name=B}] 
                \arrow[d, phantom, "\scriptstyle{\color{gray}(4)}"]
                \arrow[rr, phantom, near start, "\scriptstyle{\color{gray}(9)}"]
        & {} & {} \\
        &\alpha \beta ((\alpha \beta)^{-1} Y \otimes X)
            \arrow[d, swap, "(\eta^{\alpha\beta}_{(\alpha\beta)^{-1}Y, X})^{-1}"]
            \arrow[r, "\alpha(\eta^\beta_{(\alpha\beta)^{-1}Y, X})^{-1}"]
        & \alpha(\beta\beta^{-1}\alpha^{-1}Y \otimes \beta X)
            \arrow[d, "\alpha((\xi^\beta_{\alpha^{-1}Y})^{-1} \otimes \id)"]\\
        &\alpha \beta(\alpha \beta)^{-1} Y \otimes \alpha \beta X 
            \arrow[d, swap, "(\xi^{\alpha\beta}_Y)^{-1} \otimes \id"] 
            \arrow[ur, swap, "\eta^\alpha_{\beta(\alpha\beta)^{-1}Y, \beta X}"]
                \arrow[r, phantom, "\scriptstyle{\color{gray}(6)}"]
                \arrow[ur, bend left=10, phantom, "\scriptstyle{\color{gray}(5)}"]
        & \alpha(\alpha^{-1}Y \otimes \beta X)
            \arrow[uuu, <-, rounded corners, to path={ -- ([xshift=2ex]\tikztostart.east) 
            -| (A) node[pos=.9,right]{$\scriptstyle \alpha(\psi^\beta_{\alpha^{-1}Y})$} \tikztonodes
            |- ([xshift=2ex]\tikztotarget.east) 
            -- (\tikztotarget)}] \dar{(\eta^\alpha_{\alpha^{-1}Y, \beta X})^{-1}}\\
        &Y \otimes \alpha\beta X
            \arrow[dr, <-, swap, "\id \otimes ((\xi^\alpha_Y)^{-1} \circ \,\alpha(\xi^\beta_{\alpha^{-1}Y})^{-1} \circ \, \xi^{\alpha\beta}_Y)^{-1}"]
            \arrow[uuuuu, <-, rounded corners, to path={ -- ([xshift=-2ex]\tikztostart.west) 
            -| (C) node[pos=.9,left]{$\scriptstyle \psi^{\alpha\beta}_Y$} \tikztonodes
            |- ([xshift=-2ex]\tikztotarget.west) 
            -- (\tikztotarget)}] 
                \arrow[r, phantom, "\scriptstyle{\color{gray}(7)}"]
        & \alpha \alpha^{-1} Y \otimes \alpha \beta X
            \arrow[d, "(\xi^\alpha_Y)^{-1} \otimes \id"]
            \arrow[ul, <-, "\alpha(\xi^\beta_{\alpha^{-1}Y})^{-1} \otimes \id"] \\
        && Y \otimes \alpha\beta X
            \arrow[uuuuuuu, <-, rounded corners, to path={ -- ([xshift=2ex]\tikztostart.east) 
            -| (B) node[pos=.9,right]{$\scriptstyle (\psi^\beta_Y)^\alpha$} \tikztonodes
            |- ([xshift=2ex]\tikztotarget.east) 
            -- (\tikztotarget)}]
    \end{tikzcd}}.
            \caption{Verifying that the half-braiding behaves with composition.}
            \label{half-comp}
        \end{figure}
        Cells 1, 7, 8 and 9 commute by definition. Cells 2 and 6 commute by the distributive naturality of the tensorator. Cells 3 and 5 commute by the definition of the composition of monoidal functors. Cell 4 commutes by the functoriality of $\alpha$ and the definition of the half-braiding. Thus, the diagram commutes and we obtain:
        \begin{align*}
            & \psi_Y^{\alpha\beta} = (((\xi^{\alpha\beta})^{-1} \circ \alpha(\xi^\beta) \circ \xi^\alpha)^{-1} \circ  \otimes \id) \circ (\psi_Y^\beta)^\alpha \circ (\id \otimes ((\xi^{\alpha\beta})^{-1} \circ \alpha(\xi^\beta) \circ \xi^\alpha)) \\
            & = ((((\xi^{\alpha\beta})^{-1} \circ \alpha(\xi^\beta) \circ \xi^\alpha)^{-1} \circ ((\xi^{\alpha\beta})^{-1} \circ \alpha(\xi^\beta) \circ \xi^\alpha)) \otimes \id) \circ (\psi_Y^\beta)^\alpha \\
            & = (\id \otimes \id) \circ (\psi_Y^\beta)^\alpha = (\psi_Y^\beta)^\alpha,
        \end{align*}
        as $(\psi_Y^\beta)^\alpha$ is a half-braiding. Thus, $(\tilde\alpha, \tilde\eta^\alpha)$ is a braided monoidal autoequivalence.

        Now, we demonstrate the commutativity of the following diagram up to monoidal equivalence:
        \begin{figure}[ht]
            \begin{tikzcd}
                \underline{\Aut_\otimes}(\mc C) \rar{\underline{I^L}} \drar[swap]{\underline{I_L}} 
                & \underline{\Aut^{br}_\otimes} (\mc Z(\mc C) \vert L) \dar{\underline{F_L}} \\
                & \underline{\Aut_\otimes} (\mc Z(\mc C)_L)
            \end{tikzcd},
        \end{figure}
        
        where $\underline{F_L}$ as in \cite[Definition~4.3]{BJLP19} is strict monoidal with $\underline{F_L}(\tilde \beta, \tilde \eta^\beta, \tilde \lambda^\beta) = (\hat \beta, \hat \eta^\beta)$, $\hat\eta^\beta = \tilde \eta^\beta$, and:
        \begin{align*}
            \hat \beta(X, \psi, \mu) = \s[p]{\tilde \beta(X, \psi), \tilde\mu^\beta},
            \qquad \text{and} \qquad
            \tilde\mu^\beta := \tilde \beta(\mu) \circ \tilde \eta^\beta_{X, L} \circ (\id \otimes \tilde \lambda^\beta).
        \end{align*}
        For natural isomorphism $\tilde \pi \in \Hom((\tilde\alpha, \tilde\eta^\alpha, \tilde\lambda^\alpha), (\tilde\beta, \tilde\eta^\beta, \tilde\lambda^\beta))$, we define $\underline{F_L}(\tilde \pi)_X = \underline{F_L}(\tilde \pi_X)$.

        Define $\omega:\underline{I_L} \cong \underline{F_L} \circ \underline{I^L}$, given by:
        \begin{align*}
            & \omega_{\alpha, X}:\underline{I_L}(\alpha, \eta^\alpha)(I_L(X)) \xrightarrow{} (\underline{F_L} \circ \underline{I^L})(\alpha, \eta^\alpha)I_L(X), \\
            & \omega_{\alpha, X} = \s[p]{\bigoplus_{U \in \Irr \mc C} \eta^{\alpha}_{UX, \bar U} \circ (\eta^{\alpha}_{U, X} \otimes \id)} \circ \s[p]{\bigoplus_{U \in \Irr \mc C} \lambda^\alpha_U \otimes \id \otimes \lambda^\alpha_{\bar U} } \circ S_{I(X)},
        \end{align*}
        where $S$ is the natural isomorphism \eqref{irrperm}. It is clear that $\omega$ is a natural isomorphism forgetting the braiding and module structure. Additionally, we have:
        \begin{align*}
            \omega_{(\alpha, \eta^\alpha)} \otimes \omega_{(\beta, \eta^\beta)} = \omega_{(\alpha, \eta^\alpha) \otimes (\beta, \eta^\beta)} \circ \s[p]{\bigoplus ((\lambda^{\alpha\beta})^{-1} \circ \alpha(\lambda^\beta) \circ \lambda^\alpha) \otimes \id \otimes ((\lambda^{\alpha\beta})^{-1} \circ \alpha(\lambda^\beta) \circ \lambda^\alpha)}.
        \end{align*}
        However, $((\lambda^{\alpha\beta})^{-1} \circ \alpha(\lambda^\beta) \circ \lambda^\alpha)$ is an automorphism of a simple i.e. a scalar. By choice of component isomorphisms we obtain $\omega_{(\alpha, \eta^\alpha)} \otimes \omega_{(\beta, \eta^\beta)} = \omega_{(\alpha, \eta^\alpha) \otimes (\beta, \eta^\beta)}$.
        
        We now aim to show that $\omega$ is well defined on objects. Notice that $\alpha, \alpha^{-1}$ are biadjoint. We have:
        \begin{align*}
            & \Hom(\bar U \otimes W, \bar V)
            \cong \Hom(\alpha(\bar U \otimes W), \alpha \bar V)
            \underbrace{\cong}_{\eta^\alpha} \Hom(\alpha\bar U \otimes \alpha W, \alpha \bar V),
        \end{align*}
        meaning we have the natural identification of bases:
        \begin{align*}
            & \hat \jmath := \alpha(j) \circ \eta^\alpha, \qquad \text{and} \qquad \hat \imath^* := (\eta^\alpha)^{-1} \circ \alpha(i^*).
        \end{align*}

        Consider diagram \ref{norm-diag}.
        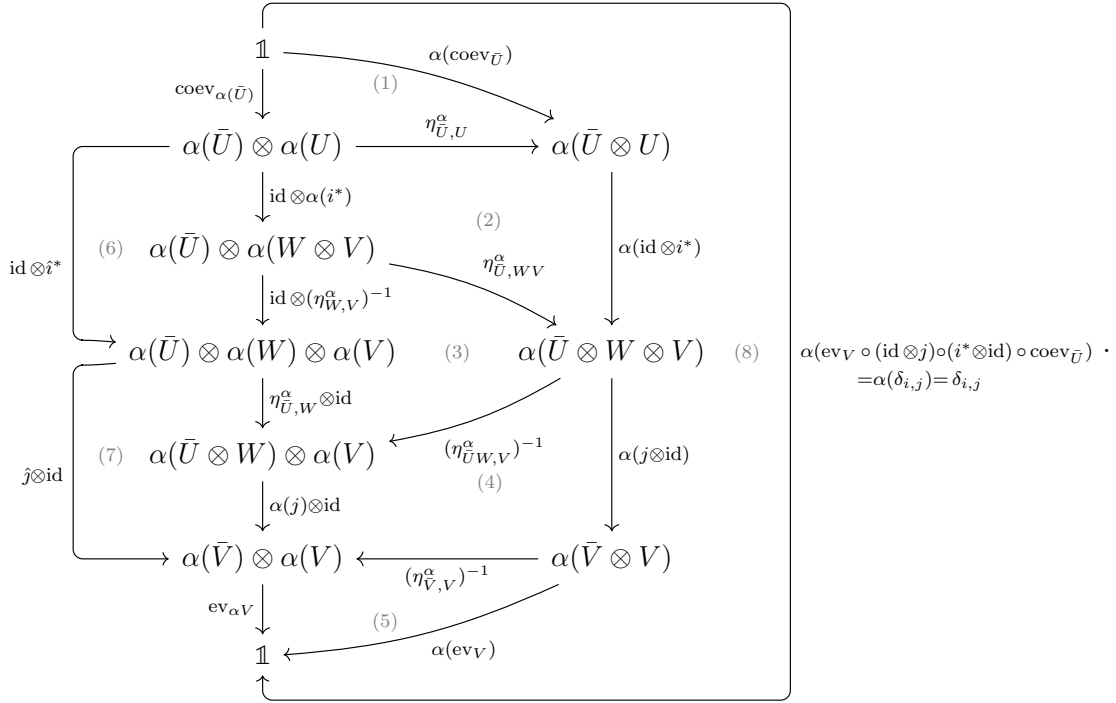
\begin{figure}[ht]
            \adjustbox{max width=0.9\textwidth}{    \begin{tikzcd}
        & [16pt] \mathbb 1 
            \arrow[d, swap, "{\coev_{\alpha(\bar U)}}"] 
            \arrow[dr, bend left=10, "{\alpha( \coev_{\bar U})}"]
                \arrow[dr, pos=0.4, phantom, "\scriptstyle{\color{gray}(1)}"]
        &[12pt] & [32pt]\\ 
        & \alpha (\bar U) \otimes \alpha (U) 
            \arrow[r, "{\eta^\alpha_{\bar U, U}} "]
            \arrow[d, "{\id \otimes \alpha(i^*)}"] 
                \arrow[ddr, bend left=20, phantom, "\scriptstyle{\color{gray}(2)}"]
        & \alpha(\bar U \otimes U) 
            \arrow[dd, "\alpha(\id \otimes i^*)"] \\
        { }
        & \alpha(\bar U) \otimes \alpha(W \otimes V) 
            \arrow[dr, bend left=10, "{\eta^\alpha_{\bar U, WV}}"]
                \arrow[l, phantom, ""{coordinate, name=B1}]
                \arrow[l, pos=0.2, phantom, "\scriptstyle{\color{gray}(6)}"]\\
        & \alpha (\bar U) \otimes \alpha (W) \otimes \alpha (V) 
            \arrow[u, swap, <-, "{\id \otimes (\eta^\alpha_{W, V})^{-1}}"] \arrow[d, "{\eta^\alpha_{\bar U, W} \otimes \id}"] 
            \arrow[uu, <-, rounded corners, to path={ -- ([xshift=-2ex,yshift=1ex]\tikztostart.west) 
            -| (B1) node[pos=.9,left]{$\scriptstyle \id \otimes \hat \imath^*$} \tikztonodes
            |- ([xshift=-2ex]\tikztotarget.west) 
            -- (\tikztotarget)}] 
                \arrow[r, phantom, "\scriptstyle{\color{gray}(3)}"]
        & \alpha(\bar U \otimes W \otimes V)
            \arrow[dd, "\alpha(j \otimes \id)"] 
                \arrow[r, phantom, ""{coordinate, name=A}]
                \arrow[r, pos=0.2, phantom, "\scriptstyle{\color{gray}(8)}"]
        & {} \\
        { }
        & \alpha(\bar U \otimes W) \otimes\alpha(V) 
            \arrow[ur, near start, <-, swap, bend right=10, "{(\eta^\alpha_{\bar U W, V})^{-1}}"]  \dar{\alpha(j) \otimes \id}
                \arrow[l, phantom, ""{coordinate, name=B2}]
                \arrow[l, pos=0.2, phantom, "\scriptstyle{\color{gray}(7)}"]\\
        & \alpha (\bar V) \otimes \alpha (V) 
            \dar[swap]{\ev_{\alpha V}} 
            \arrow[uu, <-, rounded corners, to path={ -- ([xshift=-2ex]\tikztostart.west) 
            -| (B2) node[pos=.9,left]{$\scriptstyle \hat \jmath \otimes \id$} \tikztonodes
            |- ([xshift=-2ex,yshift=-1ex]\tikztotarget.west) 
            -- (\tikztotarget)}] 
                \arrow[uur, bend right=20, phantom, "\scriptstyle{\color{gray}(4)}"]
        & \alpha(\bar V \otimes V) 
            \lar{(\eta^\alpha_{\bar V, V})^{-1}}\\
        & \mathbb 1
            \arrow[ur, bend right=10, swap, <-, "\alpha(\ev_V)"] 
            \arrow[uuuuuu, <-, swap, rounded corners, to path={ -- ([yshift=-2ex]\tikztostart.south) 
            -| (A) node[pos=1,right]{$\scriptstyle \alpha(\ev_V \circ \,(\id \otimes j) \circ (i^* \otimes \id)\, \circ\, \coev_{\bar U})$} node[pos=0.96, right] {$\qquad \scriptstyle = \alpha(\delta_{i, j}) = \, \delta_{i, j}$} \tikztonodes
            |- ([yshift=2ex]\tikztotarget.north) 
            -- (\tikztotarget)}]
                \arrow[ur, pos=0.4, phantom, "\scriptstyle{\color{gray}(5)}"]
    \end{tikzcd}}.
            \caption{A diagram verifying the normalization of a new basis.}
            \label{norm-diag}
        \end{figure}
        Cells 1 and 5 commute as $\overline{\alpha(V)} \cong \alpha(\bar V)$. Cells 2 and 4 commute by the distributive naturality of the tensorator. Cell 3 commutes by the associative naturality of the tensorator. Cells 6 and 7 commute by the definition of $\hat \imath^*, \hat \jmath$. Cell 8 commutes by the functoriality of $\alpha$ and the normalization of $i^*, j$. Thus, the diagram commutes. So, $\s{\hat \imath^*}$ and $\s{\hat \jmath}$ are bases satisfying the desired normalization.
        
        Consequently, for $\lambda^\alpha_U:Y \cong \alpha(U)$ and $\lambda^\alpha_V:Z \cong \alpha(V)$, by the choice of isomorphisms $\lambda^\alpha$ the following bases satisfy the normalization \eqref{pepsi}:
        \begin{align}
            \label{renorm}
            \s{\bar \jmath := (\lambda^\alpha_{\bar V})^{-1} \circ \alpha(j) \circ \eta^\alpha_{\bar U, \alpha^{-1}W} \circ (\lambda^\alpha_{\bar U} \otimes \xi^\alpha_W)} & \subseteq \Hom(\bar Y \otimes W, \bar Z) \\ \notag
            \s{\bar \imath^* := ((\xi^\alpha_W)^{-1} \otimes (\lambda^\alpha_V)^{-1}) \circ (\eta^\alpha_{\alpha^{-1}W, V})^{-1} \circ \alpha(i^*) \circ \lambda^\alpha_{U}} & \subseteq \Hom(Y, W \otimes Z).
        \end{align}
        
        We now aim to show that $\omega_{\alpha, X}$ is a morphism in the center. That is, that:
        \begin{align*}
            (\id \otimes \omega_{\alpha, X}) \circ \psi_{I(\alpha X), W} = \psi_{I(X), W}^\alpha \circ (\omega_{\alpha, X} \otimes \id).
        \end{align*}

        Consider diagram \ref{trans-cent}.
        \begin{figure}[ht]
            \centering
            \adjustbox{max width=0.9\textwidth}{\input{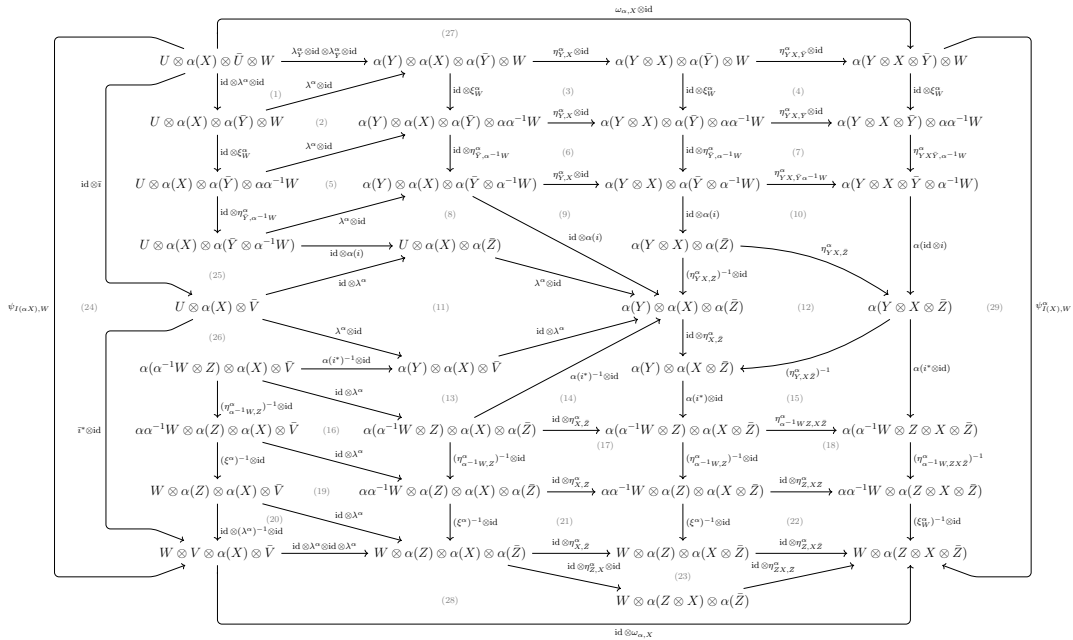}}.
            \caption{A diagram verifying that $\omega_\alpha$ is central.}
            \label{trans-cent}
        \end{figure}
        Cells 1, 2, 3, 4, 5, 6, 8, 9, 11, 13, 14, 16, 17, 19, 20, 21 and 22 commute by the naturality of the monoidal product. Cells 7, 12, 18 and 23 commute by the associative naturality of the tensorator. Cells 10 and 15 commute by the distributive naturality of the tensorator. Cells 25 and 26 commute by the definitions of $\bar \jmath, \bar \imath^*$ as in equation \eqref{renorm}. Cells 24 and 29 commute by the definition of the half-braidings $\psi_{I(\alpha X)}$, $\psi_{I(X)}^\alpha$. Cells 27 and 28 commute by the definition of the $\omega_\alpha$. Thus, the diagram commutes and $\omega_{\alpha, X}$ is a morphism in the center.
        
        We now aim to show that $\omega_{\alpha, X}$ is a morphism of $L$-modules over $\mc Z(\mc C)$. Recall from an earlier remark that $\omega_{\alpha, \mathbb1} = \tilde \lambda^\alpha$. So, we aim to show that:
        \begin{align*}
            \alpha(\eta^I_{X, \mathbb 1}) \circ \eta^\alpha_{I(X), L} \circ (\omega_{\alpha, X} \otimes \omega_{\alpha, \mathbb1}) = \omega_{\alpha, X} \circ \eta^I_{\alpha X, \mathbb 1}.
        \end{align*}
        However, by the definition of the composition of monoidal functors, $\alpha(\eta^I_{X, \mathbb 1}) \circ \eta^\alpha_{I(X), L} = \eta^{\alpha I}_{X, \mathbb 1}$. Additionally, we note that $I(\eta^\alpha_{X, \mathbb1}) = I(r_X)$ which we do not write by abuse of notation. So, if we may show that each $\omega_\alpha$ is a \emph{monoidal} natural isomorphism, we are done.

        Consider diagram \ref{nat-mon}.
        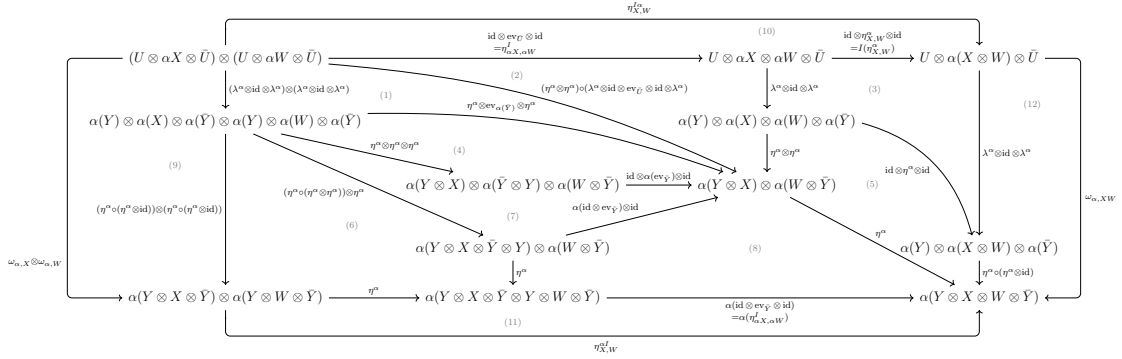
\begin{figure}[ht]
            \centering
            \adjustbox{max width=0.9\textwidth}{    \begin{tikzcd}
        & & & [12pt] {} \\ [24pt]
        & (U \otimes \alpha X \otimes \bar U) \otimes (U \otimes \alpha W \otimes \bar U)
            \arrow[rr, "\substack{\id \otimes \ev_{\bar U} \otimes \id \\ = \eta^I_{\alpha X, \alpha W}}"] 
            \arrow[d, "(\lambda^\alpha \otimes \id \otimes \lambda^\alpha) \otimes (\lambda^\alpha \otimes \id \otimes \lambda^\alpha)"] 
            \arrow[ddrr, bend left=10, "(\eta^\alpha \otimes \eta^\alpha) \circ (\lambda^\alpha \otimes \id \otimes \ev_{\bar U} \otimes \id \otimes \lambda^\alpha)"]
                \arrow[ddrr, near start, phantom, "\scriptstyle{\color{gray}(1)}"]
                \arrow[drr, bend left=5, phantom, "\scriptstyle{\color{gray}(2)}"]
        & & U \otimes \alpha X \otimes \alpha W \otimes \bar U
            \arrow[r, "\substack{\id \otimes \eta^\alpha_{X, W} \otimes \id \\ = I(\eta^\alpha_{X, W})}"] 
            \arrow[d, "\lambda^\alpha \otimes \id \otimes \lambda^\alpha"] 
                \arrow[u, phantom, ""{coordinate, name=C}]
                \arrow[u, near start, phantom, "\scriptstyle{\color{gray}(10)}"]
        & U \otimes \alpha (X \otimes W) \otimes \bar U 
            \arrow[ddd, "\lambda^\alpha \otimes \id \otimes \lambda^\alpha"] 
            \arrow[lll, <-, swap, rounded corners, to path={ -- ([yshift=2ex]\tikztostart.north) 
            |- (C) node[pos=1.3,above]{$\scriptstyle \eta^{I\alpha}_{X, W}$} \tikztonodes
            -| ([yshift=2ex]\tikztotarget.north) 
            -| (\tikztotarget)}]
                \arrow[ddd, bend left=60, near start, phantom, "\scriptstyle{\color{gray}(12)}"] \\ [12pt]
        {}
        & \alpha(Y) \otimes \alpha(X) \otimes \alpha(\bar Y) \otimes \alpha(Y) \otimes \alpha(W) \otimes \alpha(\bar Y) 
            \arrow[dr, "\eta^\alpha \otimes \eta^\alpha \otimes \eta^\alpha"] 
            \arrow[drr, bend left=10, near start, "\eta^\alpha \otimes \ev_{\alpha(\bar Y)} \otimes \eta^\alpha"] 
            \arrow[ddr, swap, "(\eta^\alpha \circ (\eta^\alpha \otimes \eta^\alpha)) \otimes \eta^\alpha"] 
            \arrow[ddd, swap, "(\eta^\alpha \circ (\eta^\alpha \otimes \id)) \otimes (\eta^\alpha \circ (\eta^\alpha \otimes \id))"]
                \arrow[l, phantom, ""{coordinate, name=A}] 
                \arrow[dr, bend left=15, near end, phantom, "\scriptstyle{\color{gray}(4)}"]
                \arrow[ddd, bend right=60, near start, phantom, "\scriptstyle{\color{gray}(9)}"]
        && \alpha(Y) \otimes \alpha(X) \otimes \alpha(W) \otimes \alpha(\bar Y) 
            \arrow[d, "\eta^\alpha \otimes \eta^\alpha"] 
            \arrow[ddr, bend left, swap, "\id \otimes \eta^\alpha \otimes \id"]
                \arrow[ur, phantom, "\scriptstyle{\color{gray}(3)}"]
                \arrow[ddr, phantom, "\scriptstyle{\color{gray}(5)}"]\\ [12pt]
        & & \alpha(Y \otimes X) \otimes \alpha(\bar Y \otimes Y) \otimes \alpha(W \otimes \bar Y)
            \arrow[r, "\id \otimes \alpha(\ev_{\bar Y}) \otimes \id"]
                \arrow[d, phantom, "\scriptstyle{\color{gray}(7)}"]
        & \alpha(Y \otimes X) \otimes \alpha(W \otimes \bar Y)
            \arrow[ddr, "\eta^\alpha"] \\ [12pt]
        & 
        & \alpha(Y \otimes X \otimes \bar Y \otimes Y) \otimes \alpha(W \otimes \bar Y)
            \arrow[ur, "\alpha(\id \otimes \ev_{\bar Y}) \otimes \id"] \arrow[d, "\eta^\alpha"]
                \arrow[dl, bend right, phantom, "\scriptstyle{\color{gray}(6)}"]
                \arrow[rr, phantom, "\scriptstyle{\color{gray}(8)}"]
        &
        & \alpha(Y) \otimes \alpha(X \otimes W) \otimes \alpha(\bar Y)
            \arrow[r, phantom, ""{coordinate, name=B}] \arrow[d, "\eta^\alpha \circ (\eta^\alpha \otimes \id)"]
        & {} \\
        & \alpha(Y \otimes X \otimes \bar Y) \otimes \alpha(Y \otimes W \otimes \bar Y)
            \arrow[r, "\eta^\alpha"] \arrow[uuuu, <-, swap, rounded corners, to path={ -- ([xshift=-2ex]\tikztostart.west) 
            -| (A) node[pos=0.6,left]{$\scriptstyle \omega_{\alpha, X} \otimes \omega_{\alpha, W}$} \tikztonodes
            |- ([xshift=-2ex]\tikztotarget.west) 
            -- (\tikztotarget)}]
        & \alpha(Y \otimes X \otimes \bar Y \otimes Y \otimes W \otimes \bar Y)
            \arrow[rr, swap, "\substack{\alpha(\id \otimes \ev_{\bar Y} \otimes \id)\\ = \alpha(\eta^I_{\alpha X, \alpha W})}"] 
                \arrow[d, phantom, ""{coordinate, name=D}]
                \arrow[d, near start, phantom, "\scriptstyle{\color{gray}(11)}"]
        &
        & \alpha(Y \otimes X \otimes W \otimes \bar Y)
            \arrow[uuuu, <-, swap, rounded corners, to path={ -- ([xshift=2ex]\tikztostart.east) 
            -| (B) node[pos=1.5,right]{$\scriptstyle \omega_{\alpha, XW}$} \tikztonodes
            |- ([xshift=2ex]\tikztotarget.east) 
            -- (\tikztotarget)}] 
            \arrow[lll, <-, swap, rounded corners, to path={ -- ([yshift=-2ex]\tikztostart.south) 
            |- (D) node[pos=0.9,below]{$\scriptstyle \eta^{\alpha I}_{X, W}$} \tikztonodes
            -| ([yshift=-2ex]\tikztotarget.south) 
            -| (\tikztotarget)}] \\ [24pt]
        & & {}
    \end{tikzcd}}.
            \caption{A diagram verifying that $\omega_\alpha$ is a monoidal natural transformation.}
            \label{nat-mon}
        \end{figure}
        Cell 1 commutes by the naturality of the monoidal product and the unitarity of $\lambda^\alpha$. Cells 2 and 3 commute by the naturality of the monoidal product. Cell 4 commutes by as $\overline{\alpha(Y)} \cong \alpha(\bar Y)$. Cells 5 and 6 commute by the associative naturality of the monoidal product. Cells 7 and 8 commute by the distributive naturality of the monoidal product. Cells 9 and 12 commute by the definition of $\omega_\alpha$. Cells 10 and 11 commute by the definition of the composition of monoidal functors. Thus, the diagram commutes and we obtain that $\omega_\alpha$ is a monoidal natural transformation. Therefore, $\omega:\underline{I_L} \xrightarrow{\sim} \underline{F_L} \circ \underline{I^L}$ is a monoidal natural isomorphism. As $\underline{I_L}$ is an equivalence, $\underline{F_L}$ is full and essentially surjective. 
        
        We now show that $\underline{F_L}$ is faithful and thus an equivalence. Given natural transformation $\pi \in \Hom(\underline{F_L}(\tilde\alpha, \tilde\eta^\alpha, \tilde\lambda^\alpha), \underline{F_L}(\tilde\beta, \tilde\eta^\beta, \tilde\lambda^\beta))$, consider diagram \ref{trans-lift}.
        
        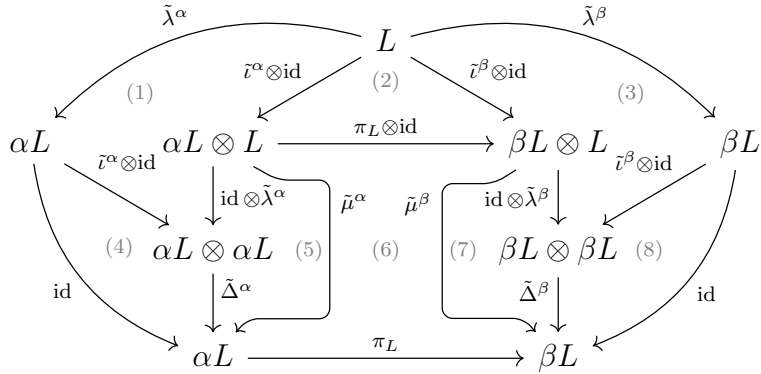
\begin{figure}[ht]
            \centering
            \adjustbox{max width=0.9\textwidth}{    \begin{tikzcd}
        && L
            \arrow[drr, bend left, "\tilde\lambda^\beta"]
            \arrow[dll, bend right, swap, "\tilde\lambda^\alpha"]
            \arrow[dl, swap, "\tilde\iota^\alpha \otimes \id"]
            \arrow[dr, "\tilde\iota^\beta \otimes \id"] 
                \arrow[dll, bend right=10, near end, phantom, "\scriptstyle{\color{gray}(1)}"]
                \arrow[d, near start, phantom, "\scriptstyle{\color{gray}(2)}"]
                \arrow[drr, bend left=10, near end, phantom, "\scriptstyle{\color{gray}(3)}"]\\
        \alpha L
            \arrow[dr, near start, "\tilde\iota^\alpha \otimes \id"]
            \arrow[ddr, bend right, swap, "\id"]
                \arrow[ddr, phantom, "\scriptstyle{\color{gray}(4)}"]
        & \alpha L \otimes L
            \arrow[d, "\id \otimes \tilde\lambda^\alpha"]
            \arrow[rr, "\pi_L \otimes \id"]
        & {} & \beta L \otimes L
            \arrow[d, swap, "\id \otimes \tilde\lambda^\beta"]
        & \beta L
            \arrow[dl, near start, swap, "\tilde\iota^\beta \otimes \id"]
            \arrow[ddl, bend left, "\id"]
                \arrow[ddl, phantom, "\scriptstyle{\color{gray}(8)}"]
            \\
        & \alpha L \otimes \alpha L
            \arrow[d, near start, "\tilde\Delta^\alpha"]
                \arrow[r, phantom, ""{coordinate, name=A}]
                \arrow[r, phantom, near start, "\scriptstyle{\color{gray}(5)}"]
                \arrow[rr, phantom, "\scriptstyle{\color{gray}(6)}"]
        & {} & \beta L \otimes \beta L
            \arrow[d, near start, swap, "\tilde\Delta^\beta"]
                \arrow[l, phantom, ""{coordinate, name=B}]
                \arrow[l, phantom, near start, "\scriptstyle{\color{gray}(7)}"]
                \\
        & \alpha L
            \arrow[rr, "\pi_L"]
            \arrow[uu, <-, swap, rounded corners, to path={ -- ([yshift=1ex]\tikztostart.north east) 
            -| (A) node[pos=1.35,right]{$\scriptstyle \tilde\mu^\alpha$} 
            \tikztonodes
            |- ([yshift=-1ex]\tikztotarget.south east)
            -- (\tikztotarget)}]
        && \beta L 
            \arrow[uu, <-, swap, rounded corners, to path={ -- ([yshift=1ex]\tikztostart.north west) 
            -| (B) node[pos=1.35,left]{$\scriptstyle \tilde\mu^\beta$} 
            \tikztonodes
            |- ([yshift=-1ex]\tikztotarget.south west)
            -- (\tikztotarget)}]
            \\
    \end{tikzcd}}
            \caption{A diagram verifying that natural transformations of modules lift.}
            \label{trans-lift}
        \end{figure}
        Cells 1 and 3 commute by the naturality of the monoidal product. Cell 2 commutes as the units reduce to $\bigoplus \alpha(\coev)$, $\bigoplus \beta(\coev)$ and as $\pi_L$ is a natural transformation. Cells 4 and 8 commute as $\alpha L$, $\beta L$ are Frobenius. Cells 5 and 7 commute by the definition of the module actions $\tilde\mu^\alpha$, $\tilde\mu^\beta$. Cell 6 commutes as $\pi_L$ is a morphism of modules. In particular, $\pi_l \circ \tilde\lambda^\alpha = \tilde \lambda^\beta$ meaning there is a canonical lift of $\pi$ to $\tilde \pi \in \Hom((\tilde\alpha, \tilde\eta^\alpha, \tilde\lambda^\alpha), (\tilde\beta, \tilde\eta^\beta, \tilde\lambda^\beta))$. $\underline{F_L}$ is thus faithful. As $\underline{I_L},$ $\underline{F_L}$ are both monoidal equivalences, then $\underline{I^L} \cong \underline{F_L}^{-1} \circ \underline{I_L}$ is an equivalence.
    \end{proof}

    \subsection{Hypergroups}

    \begin{definition}[{\cite[Definition~1.1]{SW03}}]
        A \emph{hypergroup} is a simplex $\sf{G} = \conv\s{e_0, \dots, e_{n-1}}$ equipped with an affine monoid product $*$ for which $e_0$ is the unit, together with an involutive adjunction $e_i \mapsto e_{\bar \imath}$, i.e.:
        \begin{align*}
            e_i * e_j = \sum c^k_{i, j} e_k,
            \qquad
            c^{0}_{i, j} \neq 0 \iff j = \bar \imath.
        \end{align*}
    \end{definition}

    A morphism of hypergroups is an affine morphism of monoids. We denote hypergroups using \textsf{sans serif}.

    \begin{example} 
        The convex hull of $G \hookrightarrow \mathbb C[G]$ is a hypergroup.
    \end{example}

    Let $\mc C$ a unitary fusion category. As $\mc C$ is rigid, each simple object $X$ has a unique up to isomorphism dual object $\bar X$, and $\mathbb 1$ is a subobject of $X \otimes \bar X$.

    \begin{definition}
        For $\mc C$ unitary fusion, the \emph{Grothendieck hypergroup of $\mc C$} $\sf{K_0}(\mc C)$ has basis
        \begin{align*}
            \s{\frac{[X]}{\dim X} : X \in \Irr \mc C},
        \end{align*}
        for some choice of representatives of isomorphism classes of simple objects, and monoid product given by the tensor product:
        \begin{align*}
            \frac{[X]}{\dim X} * \frac{[Y]}{\dim Y} = \frac{[X \otimes Y]}{\dim X \dim Y} = \sum_{Z \in \Irr \mc C} \frac{N^Z_{X, Y} \, \dim Z}{\dim X \dim Y} \frac{[Z]}{\dim Z}.
        \end{align*}
        The unit is $[1]$ and the adjunction sends $[X]/\dim X \mapsto [\bar X]/\dim X$.
    \end{definition}

    Let $A$ a commutative Q-system object in unitary modular category $\mc C$. One can define a \emph{convolution product} on its endomorphism space:
    \begin{align}
        \notag
        \adjustbox{scale=0.7}{    \begin{tikzpicture}[baseline=(current bounding box.center)]
        \draw (0,-0.5) node[label={east:{$A$}}] {} -- (0, 2.5);
        \node[draw, fill=white, shape=rectangle, anchor=center] at (0,1) {$f$};
        \draw (2,-0.5) node[label={east:{$A$}}] {} -- (2, 2.5);
        \node[draw, fill=white, shape=rectangle, anchor=center] at (2,1) {$g$};
        \node at (1, 1) {\Large *};
    \end{tikzpicture}} := \quad
        \adjustbox{scale=0.7}{    \begin{tikzpicture}[baseline=(current bounding box.center)]
        \draw (0,0) -- (0, 1);
        \draw (1,0) -- (1, 1);
        \draw (0, 0) arc (180:360:0.5);
        \draw (0, 1) arc (180:0:0.5);
        \draw (0.5, -0.5) node[circle, minimum size=4pt, inner sep=0pt, draw=black, fill=white] {} -- (0.5, -1) node[label={east:{$A$}}] {};
        \draw (0.5, 1.5) node[circle, minimum size=4pt, inner sep=0pt, draw=black, fill=white] {} -- (0.5, 2);
        \node[draw, fill=white, shape=rectangle, anchor=center] at (0,0.5) {$f$};
        \node[draw, fill=white, shape=rectangle, anchor=center] at (1,0.5) {$g$};
    \end{tikzpicture}} \,\,.
    \end{align}

    Further, $Q(A) := (\End(A), *, \circ)$ is a commutative semisimple algebra \cite{BD18} and admits a basis of orthogonal convolution idempotents.

    \begin{definition}
        For $A$ a commutative Q-system object in a unitary modular tensor category, define its \emph{hyperautomorphism hypergroup} $\hAut(A)$ to be the hypergroup with convex basis given by orthogonal convolution idempotents of $Q(A)$ with product given by composition of morphisms. For $\sf G$ a hypergroup, an \emph{action of $\sf G$ on $A$} is a morphism:
        \begin{align*}
            \sf G \xrightarrow{} \hAut(A).
        \end{align*}
    \end{definition}

    Other authors have studied $\hAut(A)$ under the name \emph{the symmetry hypergroup of $A$} in the context of actions on completely rational conformal nets \cite{Bis17}.

    \begin{theorem}[{\cite[Corollary~3.4]{BJ21}}]
        Let $\mc C$ a unitary fusion category, and $L = I(\mathbb 1)$ its canonical Lagrangian algebra. Then, there is an isomorphism of hypergroups:
        \begin{align}
            \notag
            \sf{K_0}(\mc C) \xrightarrow{} \hAut(L),
            \qquad \text{by} \qquad
            \frac{[X]}{\dim X}\mapsto \chi_X := \bigoplus_{U, V} \sum_i \frac{\sqrt{\dim U} \sqrt{\dim V}}{\dim X}
            \adjustbox{scale=0.7}{    \begin{tikzpicture}[baseline=(current bounding box.center)]
        \draw (0,0) node[label={below:{\strut$U$}}] {} -- (0, 2) node[label={above:{$V$}}] {};
        \draw (1,0) node[label={below:{\strut$\bar U$}}] {} -- (1, 2) node[label={above:{$\bar V$}}] {};
        \draw (0, 0.5) node[circle, minimum size=4pt, inner sep=0pt, draw=black, fill=white, label={left:{$i^\diamond$}}] {} .. controls (0, 1) and (1, 1) .. (1, 1.5) node[circle, minimum size=4pt, inner sep=0pt, draw=black, fill=white, label={right:{$i$}}] {}; 
        \draw (0.5, 1) node[label={north:{$X$}}] {};
    \end{tikzpicture}},
        \end{align}
        where $\s{i} \subseteq \Hom(X \otimes \bar U, \bar V)$, $\s{i^\diamond} \subseteq \Hom(U, V \otimes X)$ are dual with respect to the pairing:
        \begin{align} \label{copepsi}
            \delta_{i, j} = 
            \adjustbox{scale=0.7}{    \begin{tikzpicture}[baseline=(current bounding box.center)]
        \draw (0, 1.5) arc (180:0:0.5);
        \draw (0, 0.5) arc (180:360:0.5);
        \draw (0, 0.5) -- (0, 1.5);
        \draw (1, 0.5) -- (1, 1.5);
        \draw (0, 1) node[label={west:\strut$V$}] {};
        \draw (1, 1) node[label={east:\strut$\bar U$}] {};
        \draw (0.5, 1) node[label={south:$X$}] {};
        \draw (0, 0.5) node[circle, minimum size=4pt, inner sep=0pt, draw=black, fill=white, label={left:{$i^\diamond$}}] {} .. controls (0, 1) and (1, 1) .. (1, 1.5) node[circle, minimum size=4pt, inner sep=0pt, draw=black, fill=white, label={right:{$j$}}] {}; 
    \end{tikzpicture}}.
        \end{align}
    \end{theorem}

    \begin{remark} \label{canonical-action}
        The above isomorphism gives a canonical hypergroup action of $\sf{K_0}(\mc C)$ on $L$.
    \end{remark}
    
    \subsection{String Operator Symmetries}

    The author plans to elaborate the physical details of these symmetries in upcoming work.

    \begin{prop} \label{group-prop}
        The following diagram of groups is commutative:
        \begin{center}\begin{tikzcd}
            \Aut_\otimes(\mc C) \rar{I^L} \dar[swap]{\Forg} & \Aut_\otimes^{br}(\mc Z(\mc C) \vert L) \dar{\ad} \\
            \Aut (\Fus \mc C) \rar[swap]{\can} & \Aut (\End L)
        \end{tikzcd},\end{center}
        where $\ad[\tilde\alpha, \tilde\eta^\alpha, \tilde\lambda^\alpha](f) = (\tilde\lambda^\alpha)^{-1} \circ \tilde\alpha(f) \circ \tilde\lambda^\alpha$.
    \end{prop}

    \begin{proof}
        Notice that, for $\pi:(\tilde\alpha, \tilde\eta^\alpha, \tilde\lambda^\alpha) \cong (\tilde\beta, \tilde\eta^\beta, \tilde\lambda^\beta)$ two naturally isomorphism autoequivalences, we have $\tilde\lambda^\beta = \pi_L \circ \tilde\lambda^\alpha$. Thus:
        \begin{align*}
            & (\tilde\lambda^\beta)^{-1} \circ \beta(f) \circ \tilde\lambda^\beta
            = (\tilde\lambda^\alpha)^{-1} \circ \pi_L^{-1} \circ \beta(f) \circ \pi_L \circ \tilde\lambda^\alpha \\
            & = (\tilde\lambda^\alpha)^{-1} \circ \pi_L^{-1} \circ \pi_L \circ \alpha(f) \circ \tilde\lambda^\alpha
            = (\tilde\lambda^\alpha)^{-1} \circ \alpha(f) \circ \tilde\lambda^\alpha,
        \end{align*}
        so $\ad$ is indeed well defined. By \cite{BJ21}, it is sufficient to show that $\ad \circ I^L = \can \circ \Forg$ on the image of simples under the canonical action of $\Fus C$ on $L$. Let $(\alpha, \eta^\alpha) \in \underline{\Aut_\otimes}(\mc C)$. 
        
        Consider diagram \ref{group-commutes}.
        \begin{figure}[ht]
            \centering
            \adjustbox{max width=0.9\textwidth}{    \begin{tikzcd}
        &[12pt] U \otimes \bar U
            \arrow[d, "\lambda^\alpha_X \otimes \lambda^\alpha_{\bar X}"]
            \arrow[dddrr, bend left=20, "\tilde \imath^\diamond \otimes \id"]
                \arrow[ddr, phantom, "\scriptstyle{\color{gray}(3)}"]
        &[12pt]&[36pt]\\
        {} & \alpha(X) \otimes \alpha(\bar X)
            \arrow[d, "\eta^\alpha_{X, \bar X}"]
            \arrow[dr, "\alpha(i^\diamond) \otimes \id"]
                \arrow[l, phantom, ""{coordinate, name=A}]
                \arrow[l, phantom, near start, "\scriptstyle{\color{gray}(2)}"]\\
        & \alpha(X \otimes \bar X)
            \arrow[d, "\alpha(i^\diamond \otimes \id)"]
            \arrow[uu, <-, swap, rounded corners, to path={ -- ([xshift=-2ex]\tikztostart.north west) 
            -| (A) node[pos=1,left]{$\scriptstyle \tilde\lambda^\alpha$} \tikztonodes
            |- ([xshift=-2ex]\tikztotarget.west) 
            -- (\tikztotarget)}]
                \arrow[r, phantom, "\scriptstyle{\color{gray}(5)}"]
        & \alpha(Y \otimes W) \otimes \alpha(\bar X)
            \arrow[dl, "\eta^\alpha_{Y\otimes W, \bar X}"]
            \arrow[d, "(\eta^\alpha_{Y, W})^{-1} \otimes \id"] \\
        {}
        & \alpha(Y \otimes W \otimes \bar X)
            \arrow[d, "\alpha(\id \otimes i)"]
            \arrow[dr, "(\eta^\alpha_{Y, W \otimes \bar X})^{-1}"]
                \arrow[l, phantom, ""{coordinate, name=B1}]
                \arrow[l, phantom, pos=1.5, ""{coordinate, name=B2}]
                \arrow[r, phantom, "\scriptstyle{\color{gray}(6)}"]
                \arrow[l, phantom, near start, "\scriptstyle{\color{gray}(4)}"]
                \arrow[uul, phantom, bend left=50, "\scriptstyle{\color{gray}(1)}"]
        & \alpha(Y) \otimes \alpha(W) \otimes \alpha(\bar X)
            \arrow[r, "(\lambda^\alpha_Y)^{-1} \otimes \id \otimes (\lambda^\alpha_{\bar V})^{-1}"]
            \arrow[d, "\id \otimes \eta^\alpha_{W, \bar X}"]
        & U \otimes \alpha(W) \otimes V
            \arrow[dddll, bend left=20, "\id \otimes \tilde \imath"]\\
        & \alpha(Y \otimes \bar Y)
            \arrow[d, "(\eta^\alpha_{Y, \bar Y})^{-1}"]
            \arrow[uu, <-, swap, rounded corners, to path={ -- ([xshift=-2ex]\tikztostart.north west) 
            -| (B1) node[pos=1,left]{$\scriptstyle \alpha(\chi_W)$} \tikztonodes
            |- ([xshift=-2ex]\tikztotarget.south west) 
            -- (\tikztotarget)}]
                \arrow[r, phantom, "\scriptstyle{\color{gray}(7)}"]
        & \alpha(Y) \alpha(W \otimes \bar X)
            \arrow[dl, "\id \otimes i"]
                \arrow[ddl, phantom, "\scriptstyle{\color{gray}(9)}"]\\
        {} & \alpha(Y) \otimes \alpha(\bar Y)
            \arrow[d, "(\lambda^\alpha_Y)^{-1} \otimes (\lambda^\alpha_{\bar Y})^{-1}"] 
                \arrow[l, phantom, ""{coordinate, name=C}]
                \arrow[l, phantom, near start, "\scriptstyle{\color{gray}(8)}"]\\
        & V \otimes \bar V
            \arrow[uu, <-, swap, rounded corners, to path={ -- ([xshift=-2ex]\tikztostart.west) 
            -| (C) node[pos=1,left]{$\scriptstyle (\tilde\lambda^\alpha)^{-1}$} \tikztonodes
            |- ([xshift=-2ex]\tikztotarget.south west) 
            -- (\tikztotarget)}]
            \arrow[uuuuuu, <-, swap, rounded corners, to path={ -- ([xshift=-2ex]\tikztostart.south west) 
            -| (B2) node[pos=1,left]{$\scriptstyle \ad I^L[\alpha, \eta^\alpha](\chi_W)$} \tikztonodes
            |- ([xshift=-2ex]\tikztotarget.north west) 
            -- (\tikztotarget)}]
    \end{tikzcd}}.
            \caption{A diagram verifying that $\ad \circ I^L = \can \circ \Forg$.}
            \label{group-commutes}
        \end{figure}
        Cells 2 and 8 commute by the definition of $\tilde\lambda^\alpha$. Cell 4 commutes by the functoriality of $\alpha$ and the definition of $\chi_W$. Cell 1 commutes by the definition of $\ad$. Cells 5 and 7 commute by the distributive naturality of the monoidal product. Cell 6 commutes by the associative naturality of the monoidal product. We take cells 3 and 9 to be definitions of the maps $\tilde \imath^\diamond$, $\tilde \imath$. Through an analogous computation to that in \eqref{norm-diag} and \eqref{renorm}, we have that $\s{\tilde \imath^\diamond}$ and $\s{\tilde \imath}$ are bases satisfying the normalization \eqref{copepsi}. Thus, we may reindex and obtain that $\ad I^L[\alpha, \eta^\alpha](\omega_W) = \omega_{\alpha W} = \can \Forg[\alpha, \eta^\alpha](\omega_W)$ as desired.
    \end{proof}

    \begin{definition}
        Let $\mc C$ unitary fusion and $\mc D \subseteq \mc C$ a full unitary fusion subcategory. Define $\underline{\Aut_\otimes}(\mc C \vert_{\mc D})$ to be the 2-subgroup of $\underline{\Aut_\otimes}(\mc C)$ such that:
        \begin{align*}
            \underline{\Aut_\otimes}(\mc C \vert_{\mc D}) = \s{(\alpha, \eta^\alpha) \in \underline{\Aut_\otimes}(\mc C) : \alpha \vert_{\mc D} \cong \Id_{\mc D} \text{ as \emph{linear} functors}}.
        \end{align*}
    \end{definition}

    Other authors have studied the group $\Aut_\otimes(\mc C\vert_{\mc C})$ under the name \emph{soft tensor autoequivalences}, and a characterization is known for $\mc C = \mc D = \mathrm{Hilb}(G)$ \cite{Dav14}.

    \begin{definition}
        Let $A$ a commutative Q-system object in a unitary modular category $\mc C$ and $\sf {H}$ a hypergroup acting on $A$. Define $\underline{\Aut_\otimes^{br}}(\mc C \vert A, \sf H)$ to be the 2-subgroup of $\underline{\Aut_\otimes^{br}}(\mc C \vert A)$:
        \begin{align*}
            \underline{\Aut_\otimes^{br}}(\mc C \vert A, \sf H) = \s{(\alpha, \eta^\alpha, \lambda^\alpha) \in \underline{\Aut_\otimes^{br}}(\mc C \vert A) : \lambda^\alpha \circ e_i = \alpha(e_i) \circ \lambda^\alpha, \text{ for all } e_i \in \sf{H}}.
        \end{align*}
    \end{definition}

    \begin{cor}
        For $\mc C$ a unitary fusion category and $\mc D \subseteq \mc C$ a full fusion subcategory, there is an equivalence of 2-groups:
        \begin{align*}
            \underline{I^L}\vert_{\mc D}:\underline{\Aut_\otimes}(\mc C\vert_{\mc D}) \xrightarrow{\sim} \underline{\Aut_\otimes^{br}}(\mc Z(\mc C) \vert L, \Fus \mc D)
        \end{align*}
    \end{cor}

    \begin{proof}
        We remark that there is a well defined notion of $\Fus \mc D$ acting on $L$ by restricting the canonical action. Consider the subset of automorphisms of $\End(L)$ fixing $\s{\chi_W : W \in \Irr \mc D}$. The preimage of these automorphisms under $\ad$ will be $[\tilde\alpha, \tilde\eta^\alpha, \tilde\lambda^\alpha]$ such that:
        \begin{align*}
            \ad[\tilde\alpha, \tilde\eta^\alpha, \tilde\lambda^\alpha](\chi_W) = \chi_W
            \iff \tilde\alpha(\chi_W) \circ \tilde\lambda^\alpha = \tilde\lambda^\alpha \circ \chi_W,
        \end{align*}
        which is exactly $\Aut_\otimes^{br}(\mc Z(\mc C) \vert L, \Fus \mc D)$. Similarly, the preimage of such automorphisms under $\can \circ \Forg$ are the autoequivalences fixing simples in $\mc D$, which is exactly $\Aut_\otimes(\mc C\vert_{\mc D})$. By proposition \ref{group-prop}, $I^L:\Aut_\otimes(\mc C\vert_{\mc D}) \xrightarrow{\sim} \Aut_\otimes^{br}(\mc Z(\mc C) \vert L, \Fus \mc D)$ are isomorphic as groups. As 2-subgroups are `full,' we obtain the categorified result.
    \end{proof}

    \subsubsection{Examples}

    While a general description of the canonical Lagrangian algebra and its braiding are given, a decomposition into simple objects is often difficult. In the following examples we do not give the associator for the automorphism 2-group, and will only describe the simple objects of $\mc Z(\mc C)$ when a general description is known.

    \begin{example}
        Let $G$ a finite group and consider the category $\mathrm{Hilb}(G)$ of finite dimensional $G$-graded Hilbert spaces. Full fusion subcategories of $\mathrm{Hilb}(G)$ are in one-to-one correspondence with subgroups $H$ of $G$. It is well known \cite[Proposition~2.6.1]{EGNO} that:
        \begin{align*}
            \Aut_\otimes(\mathrm{Hilb}(G)) \cong \Aut(G) \rtimes H^2(G, \mathbb C^\times).
        \end{align*}

        It follows that:
        \begin{align*}
            \Aut_\otimes(\mathrm{Hilb}(G) \vert_{\Vec(H)}) \cong \Stab(N) \rtimes H^2(G, \mathbb C^\times).
        \end{align*}

        The simple objects of $\mc Z(\mathrm{Hilb}(G))$ are indexed by conjugacy classes together with representations of centralizers of orbits \cite[Example~8.5.4]{EGNO}.

        Let $G = S_4$. It is known that every automorphism of $S_4$ is inner, so stabilizers under the inner action are exactly centralizers. Additionally, $H^2(S_4, \mathbb C^\times) \cong \mathbb Z_2$ has no nontrivial automorphism. There are 30 subgroups of $S_4$, of which there are 11 distinct types. A representative and the number of each type is given.
        \begin{itemize}
            \item The trivial subgroup, $\s{0}$ (1).
            \item $\mathbb Z_2$ generated by a single transposition, $\s[g]{(12)}$ (6).
            \item $\mathbb Z_2$ generated by a pair of disjoint transpositions, $\s[g]{(12)(34)}$ (3).
            \item $\mathbb Z_3 \cong A_3$ generated by a 3-cycle, $\s[g]{(123)}$ (4).
            \item $\mathbb Z_4$ generated by a 4-cycle, $\s[g]{(1234)}$ (3). 
            \item Normal Klein four group $V_4 := \s{0, (12)(34), (13)(24), (14)(23)}$ (1).
            \item Non-normal Klein four groups, $\s[g]{(12), (34)}$ (3).
            \item Dihedral group $D_8$, $\s[g]{(1234), (13)}$ (3).
            \item Symmetric group $S_3$, $\s[g]{(12), (23)}$ (4).
            \item Alternating group $A_4$ (1).
            \item Symmetric group $S_4$ (1).
        \end{itemize}

        One may explicitly compute the stabilizers:
        \begin{itemize}
            \item $\Stab(0) = S_4$
            \item $\Stab((12)) = \Stab((12), (34)) = \s[g]{(12), (34)}$
            \item $\Stab((12)(34)) = \s[g]{(1324), (12)}$
            \item $\Stab((123)) = \s[g]{(123)}$
            \item $\Stab((1234)) = \s[g]{(1234)}$
            \item $\Stab(V_4) = V_4$
            \item $\Stab((1234), (13)) = \s[g]{(13)(24)}$
            \item $\Stab(S_3) = \Stab(A_4) = \Stab(S_4) = 0$
        \end{itemize}
    \end{example}

    \begin{example}
        Let $A$ a finite Abelian group, $\chi$ a nondegenerate bicharacter on $A$, and $\tau$ a choice of sign. Denote the associated Tambara-Yamagami category \cite{TY98} by $\mc{TY}(A, \chi, \tau)$. It is well known \cite[Lemma~2.16]{Edi22} that the group of isomorphism classes of monoidal autoequivalences of $\mc{TY}(A, \chi, \tau)$ is isomorphic to the group of automorphisms of $A$ preserving the bicharacter $\chi$, i.e.:
        \begin{align*}
            \Aut_\otimes(\mc{TY}(A, \chi, \tau)) \cong \Aut(A, \chi).
        \end{align*}

        Further, it is known that $\Aut(\Id_{\mc C}) = \Hom(\mc U_{\mc C}, \mathbb C^\times)$ for $\mc U_{\mc C}$ the universal grading group \cite[Proposition~4.14.3]{EGNO} and that $\mc U_{\mc C} = \mathbb Z_2$ \cite[Proposition~5.3]{Nat13}. The associator for $\underline{\Aut_\otimes}(\mc C)$ is given by a 2-cocycle $\omega \in H^2(\Aut(A, \chi), \mathbb Z_2)$. 
        
        The strict full fusion subcategories of $\mc C$ are in one-to-one correspondence with subgroups $H \subseteq A$, which we denote $\s[g]{H} \subseteq \mc C$. It follows that:
        \begin{align*}
            \Aut_\otimes(\mc C \vert_{\s[g]{H}}) \cong \Stab(H, \chi) \subseteq \Aut(A, \chi),
        \end{align*}
        the subgroup of autoequivalences preserving the bicharacter $\chi$ that fix the subgroup $H$.

        A description of the simple objects of $\mc Z(\mc{TY}(A, \chi, \tau))$ is known \cite[Section~4]{GNN09}, of which there $2\s[v]{A}$ invertible simple objects and $4\s[v]{A} + \binom{\s[v]{A}}2$ total simples. The author of this note was unable to find a description of the simple objects of the center of arbitrary $\mathcal{TY}(A, \chi, \tau)$ nor of its Brauer-Picard group.
    \end{example}

    \begin{example}
        Consider the category $\mc C(\mathfrak{so}_{2r+1}, 2)$ of $\mc U_q(\mathfrak{so}_{2r+1})$ tilting modules at level $k = 2$ \cite{Sch20}. This category has rank $4+r$ and Frobenius-Perron dimension $4 + 8r$. The simple roots and fundamental weights of $\mathfrak{so}_{2r+1}$ are: 
        \begin{align*}
            (\text{simple roots}) & = \s{\alpha_i = \epsilon_i - \epsilon_{i+1} : i < r} \cup \s{\alpha_r = \epsilon_r} \\
            (\text{fundamental weights}) & = \s{\lambda_i = \sum_{j \leq i} \epsilon_j : i < r} \cup \s{ \lambda_r = \frac12 \sum_{j \leq r}\epsilon_j}
        \end{align*}
        
        The unique long root in the Weyl chamber is $\theta := \alpha_1 + 2\alpha_2 + \dots + 2 \alpha_r$. The simple objects of $\mc C(\mathfrak{so}_{2r+1}, 2)$ correspond to weights in the Weyl alcove, which consists of weights $\lambda$ in the Weyl chamber for which $\s[g]{\lambda, \theta} \leq k = 2$. Explicitly, these are the weights:
        \begin{align*}
            \Lambda_0 = \s{0, \lambda_1, \lambda_2, \dots, \lambda_r, 2\lambda_1, 2\lambda_r, \lambda_1 + \lambda_r}.
        \end{align*}

        The simply connected Lie group corresponding to $\mathfrak{so}_{2r + 1}$ is the spin group ${\mathrm{Spin}(2r + 1)}$, whose center is $\mathbb Z_2$ \cite[Section~5.3]{Var04}. Full fusion subcategories of $\mc C(\mathfrak{so}_{2r+1}, 2)$ correspond to closed subsets of the Weyl alcove $\Lambda_0$. There are 4 closed subsets geometrically arising from subgroups of $Z(\mathrm{Spin}(2r+1))$ and possibly many anomalous subsets arising from divisors of $2r+1$ \cite{Saw06}. Explicitely these are:
        \begin{align*}
            & \Delta_0 = \s{0}, \qquad
            \Delta_{\mathbb Z_2} = \s{0, 2\lambda_1}, \qquad
            \Gamma_0 = \Lambda_0, \\
            & \Gamma_{\mathbb Z_2} = \s{0, \lambda_1, \lambda_2, \dots, \lambda_{r-1}, 2\lambda_1, 2\lambda_r} = \Lambda_0 \setminus \s{\lambda_r, \lambda_1 + \lambda_r}, \\
            & \Xi_j := \s{0, 2\lambda_1, \lambda_j, \lambda_{2j}, \dots, \lambda_{(2r+1-j)/2}}, \quad \text{ for } \quad 2 < j \vert 2r+1.
        \end{align*}

        Define $\omega(n)$ to be the number of distinct prime factors of $n$. The group structure of the autoequivalences of $\mc C(\mathfrak{so}_{2r+1}, 2)$ is known \cite{Edi22}:
        \begin{align*}
            \Aut_\otimes(\mc C(\mathfrak{so}_{2r+1}, 2)) \cong
            \begin{cases}
                \mathbb Z_2 \times \mathbb Z_2 \times \mathbb Z_2^{\omega(2r + 1) - 1}
                & p \equiv 1 \mod 4 \text{ for all prime } p \vert 2r+1\\
                \mathbb Z_2 \times \mathbb Z_2^{\omega(2r + 1) - 1}
                & \text{otherwise}.
            \end{cases}
        \end{align*}
        The first $\mathbb Z_2$ corresponds to the autoequivalence which swaps $\lambda_r$ and $\lambda_1 + \lambda_r$ but fixes every other simple. The second and third factor corresponds to exotic fusion ring automorphisms, which act by sending:
        \begin{align*}
            \lambda_i \mapsto \lambda_{\min(mi \mod 2r+1, -mi \mod 2r+1)}, \quad \text{ for } \quad m \in \mathbb Z_{2r+1}^\times, \quad m^2 \equiv \pm 1 \mod 2r+1.
        \end{align*}

        Write $2r+1 = \prod p_i^{\alpha_i}$ a decomposition into distinct prime powers. In particular:
        \begin{align*}
            n^2 \equiv 1 \mod 2r+1 \iff 
            n^2 \equiv 1 \mod p_i^{\alpha_i}, \,\, \forall i \iff
            n \equiv \pm 1 \mod p_i^{\alpha_i}, \,\, \forall i.
        \end{align*}
        
        There are $2^{\omega(2r+1)-1}$-many choices of $n$ modulo sign. Let $j$ a divisor of $2r+1$ and define:
        \begin{align*}
            \theta(j;2r+1) := \#\s{i \in \omega(2r+1) : p_i^{\alpha_i} \text{ divides } j}.
        \end{align*}

        One sees that $\Xi_j$ is fixed by automorphisms for which ${n \equiv 1 \mod (2r+1 / j)}$, of which there are $2^{\theta(j;2r+1)}$-many possible choices modulo sign. It follows that:

        \begin{itemize}
            \item $\Stab(\Delta_0) = \Stab(\Delta_{\mathbb Z_2}) = \Aut_\otimes(\mc C(\mathfrak{so}_{2r+1}, 2))$,
            \item $\Stab(\Gamma_0) = 0$,
            \item $\Stab(\Gamma_{\mathbb Z_2}) = \mathbb Z_2$,
            \item $\Stab(\Xi_j) = \mathbb Z_2 \times \mathbb Z_2^{\theta(j;2r+1)}$.
        \end{itemize}
    \end{example}

    \printbibliography

\end{document}